%% file: nprk-intro-paper-shortened.tex
\documentclass[final]{siamart171218}
\pdfoutput=1

\usepackage{amsmath,amssymb,amsfonts,mathtools}
\usepackage{comment}
\usepackage{enumerate}
\usepackage{multirow}
\usepackage{hyperref}
\usepackage{cleveref}
\hypersetup{colorlinks=true, urlcolor=blue,}
\usepackage{xcolor}
\usepackage{caption}
\usepackage{subcaption}
\usepackage{graphbox}
\usepackage[title]{appendix}
\newtheorem{remark}{Remark}
\usepackage{dashrule}
\usepackage{enumitem}

\definecolor{plot_black}{rgb}{0,0,0}
\definecolor{plot_red}{rgb}{0.7529, 0.2235, 0.1686}
\definecolor{plot_blue}{rgb}{0.1608, 0.5020, 0.7255}

\usepackage{algorithm,algorithmicx}
\usepackage[noend]{algpseudocode}
\newcommand*\Let[2]{\State #1 $\gets$ #2}
\newcommand{\pluseq}{\mathrel{+}=}

\title{A New Class of Runge-Kutta Methods for Nonlinearly Partitioned Systems\thanks{Corresponding author email address: tbuvoli@tulane.edu
\funding{B. S. S. was supported by the Laboratory Directed Research and Development program of Los Alamos National Laboratory as a Nicholas C. Metropolis Fellow and under project number 20220174ER. Los Alamos National Laboratory report number LA-UR-24-20116. T. B. was funded in part by the NSF grant OIA-2327484.}}}

\author{Tommaso Buvoli\thanks{Mathematics Department, Tulane University, New Orleans, LA 70118, USA}
   \and Ben S. Southworth\thanks{Theoretical Division, Los Alamos National Laboratory, USA.}
}

\begin{document}
\maketitle

\begin{abstract}
This work introduces a new class of Runge-Kutta methods for solving nonlinearly partitioned initial value problems. These new methods, named nonlinearly partitioned Runge-Kutta (NPRK), 
generalize existing additive and component-partitioned Runge-Kutta methods, and allow one to distribute different types of implicitness within nonlinear terms. The paper introduces the NPRK framework and discusses order conditions, linear stability, and the derivation of implicit-explicit and implicit-implicit NPRK integrators. 
The paper concludes with numerical experiments that demonstrate the utility of NPRK methods for solving viscous Burger's and the gray thermal radiation transport equations.
\end{abstract}

% ------------------------------------------------------- %
% ------------------------------------------------------- %

\input{section-introduction.tex}

% ------------------------------------------------------- %
% ------------------------------------------------------- %

\input{section-background.tex}

% ------------------------------------------------------- %
% ------------------------------------------------------- %

\input{section-motivation.tex}

% ------------------------------------------------------- %
% ------------------------------------------------------- %

\input{section-nprk.tex}

% ------------------------------------------------------- %
% ------------------------------------------------------- %

\input{section-example-methods.tex}

% ------------------------------------------------------- %
% ------------------------------------------------------- %

\input{section-algorithm-imex.tex}

% ------------------------------------------------------- %
% ------------------------------------------------------- %

\input{section-numerical-experiments.tex}

% ------------------------------------------------------- %
% ------------------------------------------------------- %

\input{section-conclusion.tex}

% ------------------------------------------------------- %
% ------------------------------------------------------- %

\section*{Acknowledgements}  
We would like to thank Ryosuke Park and Marc Charest as the primary developers of the \emph{flecsim} package we used for thermal radiation transport calculations. This research used resources provided by the Darwin testbed at Los Alamos National Laboratory (LANL) which is funded by the Computational Systems and Software Environments subprogram of LANL's Advanced Simulation and Computing program (NNSA/DOE). Lastly, T. Buvoli was funded by the National Science Foundation grant NSF-OIA-2327484.

\input{section-appendix.tex}

\bibliographystyle{plain}
\bibliography{refs}{}

\newpage
\input{section-supplemental.tex}

\end{document}

%% file: section-introduction.tex
\section{Introduction}\label{sec:intro}

We are interested in the efficient and accurate numerical solution of
ordinary differential equations (ODEs) arising from the spatial
discretization of partial differential equations (PDEs) and partial
integro-differential equations (PIDEs). Such problems typically give rise to a combination of numerical challenges as a result of multiple scales, multiple
stiffnesses, and nonlinear coupling between variables and/or scales
of the problem.

The solution of these equations amounts to solving an initial value problem (IVP)
\begin{align}
	y' = G(y), \quad y(t_0) = y_0.
	\label{eq:generic-ivp}
\end{align}
For the equations we are interested in, explicit time integration is computationally intractable due to stability constraints arising from stiff processes. A classical approach is to apply an implicit time integration method to the full right-hand-side $G(y)$. Although this strategy is highly effective at overcoming stability limitations, the development of efficient solvers presents its own set of challenges (e.g., see \cite{Elman.2014,Vassilevski.2008}). 
An alternative idea is to only treat the stiff components of $G(y)$ implicitly. Suppose, for example that one can separate the stiff and non-stiff components of $G(y)$ either by separating additively combined terms in $G(y)$ or by partitioning the components of $x$. In either case, one obtains the equation
\begin{align}
	y' = G^{\{1\}}(y) + 	G^{\{2\}}(y)
	\label{eq:additive-ivp}
\end{align}
(see \cref{subsubsec:additive-component-equivalence}). For such additively partitioned equations, it can be advantageous to only treat the first component $G^{\{1\}}(y)$ implicitly (e.g. consider the canonical diffusion-reaction system where only diffusion is treated implicitly). This leads to the well-known class of implicit-explicit (IMEX) integrators \cite{Ascher.1995,Ascher.1997} and the more general class of additive integrators that allow for $n$ partitions and more varied forms of implicitness \cite{Kennedy.2003tv4,Kennedy.2019,Sandu.2015}.

This paper is concerned with systems of ODEs that cannot be additively or component-wise partitioned. To capture this generality we re-express the function $G(y) : \mathbb{R}^n \to \mathbb{R}^n$ as a multivalued function $F(u,v) : \mathbb{R}^n \times \mathbb{R}^n \to \mathbb{R}^n$ that satisfies the equality $G(y) = F(y,y)$; a simple scalar example is $G(y)=y^2$ and $F(u,v) = uv$. The primary aim of this work is to introduce a new class of
nonlinearly partitioned Runge-Kutta (NPRK) methods for solving the IVP \cref{eq:generic-ivp} rewritten as
\begin{align}
	y' = F(y,y),
	\label{eq:np-ivp}
\end{align}
This approach makes it possible to treat each argument of $F$ with a differing degree of implicitness, and leads to nonlinear generalizations of existing additive implicit-explicit (IMEX), explicit-explicit (EXEX), and implicit-implicit (IMIM) integrators. Moreover, the proposed NPRK framework can be seen as a
generalization of the pioneering semi-implicit work of Boscarino et al. \cite{Boscarino.2016,Boscarino.2015}, and recently \cite{Sebastiano.2023}, where the governing equations are copied to facilitate
implicit-explicit-like integration for nonlinearly coupled systems. 

The nonlinearly partitioned ODE \cref{eq:np-ivp} has been successfully applied  to complex flow problems in a number of recent works, e.g., \cite{Boscarino.2022,Boscheri.2022,southworth2023implicit}. In this work we will demonstrate how our proposed NPRK methods can be used to solve the thermal radiative transfer (TRT) equation, a complex and high-dimensional PIDE that governs the evolution of radiation energy and temperature, and is critical to complex physics such as inertial confinement fusion \cite{Haines.2022}. In non-relativistic regimes, the timescales of interest are typically much larger than explicit stability constraints. However, the high dimensionality of the underlying PIDE and its various types of stiffness render classical implicit integration extremely expensive \cite{Park.2013}. An additional challenge arises from certain reaction or ``opacity'' terms that depend nonlinearly on temperature and require tabular data lookups for realistically representing high-energy regimes; incorporating these terms in a nonlinear implicit solve is nontrivial.
Furthermore, in recent work we developed a nonlinear partitioning of the TRT equations to greatly reduce the computational cost of time integration by only implicitly integrating certain stiff asymptotic limits \cite{imex-trt}. Neither of these formulations are amenable to standard additive partitioning, but can be naturally represented as a nonlinear partition 
\cref{eq:np-ivp}. In future work, we will also apply the proposed NPRK methods to the coupling of radiation with hydrodynamics \cite{Lowrie.1999} using a nonlinear partitioning developed in our recent work \cite{southworth2023implicit}. This nonlinear partitioning facilitates high-order time integration compared with a previously used first-order Lie-Trotter-like operator split that was specially designed for the radiation hydrodynamics equations, e.g., see \cite{Bates.2001,southworth2023implicit} (Lie-Trotter-like in the sense that it is not a classical additive Lie-Trotter operator split, and various nonlinearly coupled variables are ``frozen'' at subsequent points in the algorithm).

This paper proceeds as follows.  
\Cref{sec:background} provides background on existing partitioned time integration methods and \Cref{sec:motivation} motivates nonlinear partitioning using the viscous Burger's equation. 
Our new class of NPRK methods is then introduced in \Cref{sec:gsi}, including the relationship with existing Runge-Kutta methods, order conditions, and  linear stability analysis. Complete order conditions via rooted tree analysis is the subject of a companion paper \cite{nprk2}. Example methods of orders $1-3$ are derived in \Cref{sec:example}, and a practical algorithm for a subclass of NPRK methods is detailed in \Cref{sec:alg}. We then demonstrate the new methods in \Cref{sec:results}, first on a model viscous Burger's equation, and then on a nonlinear partitioning of gray TRT. Conclusions and future work are discussed in \Cref{sec:conclusions}.

%% file: section-background.tex
\section{Background and related integration techniques}
\label{sec:background}

We begin by discussing existing integration techniques for solving partitioned equations, including integrators for the nonlinearly partitioned ODE \cref{eq:np-ivp}. We will revisit these methods in later sections to highlight connections with the newly proposed NPRK methods. 

\subsection{Additive and partitioned integrators}
	
	Two commonly used types of partitioning are additive and component-wise partitioning. 
	We provide a brief overview of these concepts and highlight their equivalence.
	
	\subsubsection{Additive partitioning}
	
	When the right-hand-side of \cref{eq:generic-ivp} can be expressed as the sum of multiple functions, then it is an additively partitioned system. A two-part additively partitioned system is
	\begin{align}
    		y' &= G^{\{1\}}(y) + G^{\{2\}}(y) &
    		y(t_0) &= y_0
    	\label{eq:additive}
	\end{align}
	It is natural to consider integrating each component with a different Runge-Kutta method. This leads to the additive Runge-Kutta (ARK) methods that were first introduced in \cite{Cooper.1983sej} and later popularized in \cite{Ascher.1997,Kennedy.2003tv4}. An $s$-stage ARK method is
	\begin{align}
   		\begin{aligned}
   		\text{stage values:} && Y_{i} &= y_n + h \sum_{j=1}^s \left[ a^{\{1\}}_{ij} G^{\{1\}}(Y_j) + a^{\{2\}}_{ij} G^{\{2\}}(Y_j) \right], \quad i=1,\ldots, s\\
   		\text{output:} && y_{n+1} &= y_n + h \sum_{i=1}^s \left[ b^{\{1\}}_{i} G^{\{1\}}(Y_i) + b^{\{2\}}_{i} G^{\{2\}}(Y_i) \right].
  		\end{aligned}
  		\label{eq:nprk-underlying-ark}
	\end{align}
	If the matrix $a_{ij}^{\{1\}}$ is lower triangular (with non-zero diagonals), and $a_{ij}^{\{2\}}$ is strictly lower triangular then the ARK method is an implicit-explicit (IMEX) method that treats the first component implicitly and the second explicitly. Many papers have explored ARK methods (typically of the IMEX
variety) \cite{Pareschi.2000,Kanevsky.2007}, focusing on topics such
as increased order and stability \cite{Izzo.2017,Kennedy.2019, Calvo.2001},
strong stability preservation \cite{Conde.2017,Higueras.2014,Higueras.2006},
reduced memory \cite{Cavaglieri.2015,Higueras.2016}, and hyperbolic
equations with stiff relaxation \cite{Pareschi.2005,Boscarino.2013,Dimarco.2013}. The order conditions for ARK coefficients are also well-known; see for example \cite{Kennedy.2003tv4}.

	\subsubsection{Component partitioning}
	
	When the solution vector can be naturally partitioned into separate components, one obtains a component partitioned equation. For example, a two-component partitioned equation is 
	\begin{align}
       	z' &= H^{\{1\}}(z,w), \hspace{1ex} w' = H^{\{2\}}(z,w), &&
    	z(t_0) = z_0, \hspace{1ex} w(t_0) = w_0
    	\label{eq:partitioned}
	\end{align}
	which is the natural form for canonical Hamiltonian systems. Partitioned Runge-Kutta (PRK) methods for solving \cref{eq:partitioned} were first derived in \cite{Rice.1960,hofer1976partially}
(see also \cite[II.15]{Hairer.1993}). PRK methods integrate each
component with a different Runge-Kutta method and have the ansatz
\begin{align}
    \begin{small}
    \begin{aligned}
        \text{stages:} \quad &  \left\{ 
        \begin{aligned}
            Z_{i} &= z_{n} + h \sum_{j=1}^s a_{ij} H^{\{1\}}(Z_{j}, W_{j}), \\
            W_{i} &= w_{n} + h\sum_{j=1}^s \widehat{a}_{ij} H^{\{2\}}(Z_{j}, W_{j})  \\ 
        \end{aligned} \right.
        \hspace{5em} i = 1, \ldots, s, \\
        \text{output:} \quad & \left\{  
         \begin{aligned}
            z_{n+1} &= z_{n} + h \sum_{j=1}^s b_{j} H^{\{1\}}(Z_{j}, Y_{j}), &&
            w_{n+1} = w_{n} + h\sum_{j=1}^s \widehat{b}_{j} H^{\{2\}}(Z_{j}, W_{j}).\\
        \end{aligned} \right.
    \end{aligned}
    \end{small}
    \label{eq:partitioned-runge-kutta}
\end{align}
The order conditions for PRK methods are known and can be found in \cite{Hairer.1993,Jackiewicz.2000}.
	
	\subsubsection{Equivalence of additive and partitioned integration}
	\label{subsubsec:additive-component-equivalence}
		
	Although additive and component partitioning appear different, the
	formulations are equivalent, in the sense that either can be expressed
	as the other. For example, if we let $y = [z; w]$,
	$G^{\{1\}}(y) \coloneqq [H^{\{1\}}(z,w); 0]$, and $G^{\{2\}}(y) \coloneqq [0; H^{\{2\}}(z,w)]$,
	then we can write \cref{eq:partitioned} as \cref{eq:additive}. Going the other way, for any $z$ and $w$ satisfying $y = z+w$, we define
	$H^{\{1\}}(z,w) = G^{\{1\}}(z+w)$, and
	$H^{\{2\}}(z,w) = G^{\{2\}}(z+w)$ then \cref{eq:additive} is equivalent to
	\cref{eq:partitioned}. Note that 
	in this latter reformulation, we have doubled the number of
	unknowns, and the initial conditions are no longer unique ($w_0$ can take any value so long as $z_0 = y_0 - w_0$).
	\subsubsection{Generalizations}
	
	Additive and component-partitioned systems can be generalized to an arbitrary
number of pieces; for additive this implies that $G(y)$ is the sum of $m$ functions, while for component splittings $y$ is divided into $m$ components. This leads to the general class of ARK methods for an arbitrary number of partitions \cite{Kennedy.2003tv4}. More recently, generalized Additive Runge-Kutta (GARK) methods were proposed in \cite{Sandu.2015}, which generalize classical ARK methods by allowing for
different stage values as arguments of different components of the
right-hand side \cite{Sandu.2015}. GARK methods can be
written as ARK methods, but the resulting ARK-form is cumbersome.
The GARK framework has led to the development of many new methods,
e.g., \cite{Gunther.2016,Sandu.2021,Sandu.2019,Gunther.2021},
particularly multirate integrators. However, despite
the added flexibility of GARK methods over ARK, it does not address
the objectives of this paper, namely time integrators that are naturally amenable
to nonlinear partitioned equations.

\subsection{Semi-implicit RK}
\label{subsec:sirk}

Semi-implicit RK (SIRK) methods (also referred to as linearly implicit
IMEX RK methods) were developed in \cite{Boscarino.2016,Boscarino.2015}, and are to
our knowledge the only general class of methods in the literature
that can be directly applied to the nonlinearly partitioned equation \cref{eq:np-ivp}. The SIRK framework proposed in \cite{Boscarino.2016,Boscarino.2015} generalizes the initial semi-implicit RK methods from \cite{Zhong.1996} by introducing a novel way to apply IMEX-ARK methods. 
Specifically, they propose to rewrite equation \cref{eq:generic-ivp} as  \cref{eq:np-ivp} and assume that $F$ has a stiff dependence
on the first variable and a nonstiff dependence on the second variable.
An artificial partitioned system is then posed as 
\begin{equation}\label{eq:partitioned-simex}
    {y^*}' = F(y,y^*), \hspace{3ex} y' = F(y,y^*),
\end{equation}
where the continuous solutions $y(t) = y^*(t)$, but we treat the
partitioning explicitly in $y^*$ and implicitly in $y$. Note that we
now have a component partitioned set of equations as a specific case
of \cref{eq:partitioned}, and as a result can apply standard PRK 
methods \cref{eq:partitioned-runge-kutta} to \cref{eq:partitioned-simex}.

A general SIRK method for solving \cref{eq:partitioned-simex}  is a PRK method \cref{eq:partitioned-runge-kutta} defined by a double Butcher tableaux with equal RK weights $\mathbf{b} = \widehat{\mathbf{b}}$, upper triangular $A^{\{1\}} = A^{I}$ and lower triangular $A^{\{2\}}=A^{E}$.
The resulting method treats $y^*$ explicitly using $\{c^E, A^E, b\}$ and $y$ implicitly using $\{c^I, A^I, b\}$. 

\subsection{Other notable integration methods}

Rosenbrock methods \cite{rosenbrock1963some,wanner1977integration,steihaug1979attempt,
Podhaisky.2005,Lang.2021} are implicit methods for solving \cref{eq:generic-ivp} that avoid nonlinear solves in favor of inverting a single linear system at each stage. For classical Rosenbrock methods this system involves the exact Jacobian of $G(y)$, while Rosenbrock-W methods allow for an inexact Jacobian. Although Rosenbrock methods can offer improved computational efficiency compared to implicit RK methods, they do not facilitate flexible explicit, implicit, and linearly implicit treatment of variables within nonlinear terms. 
In particular, a Rosenbrock method imposes the same (linearly implicit) approach uniformly across a problem.

Explicit-Implicit-Null (EIN) methods \cite{duchemin2014explicit,wang2020local,tan2022stability} add and subtract an arbitrary operator to the equation right-hand-side, then treat the added operator implicitly and all remaining terms explicitly. This can be interpreted as an IMEX-ARK method applied to the additive splitting \cref{eq:additive} with $G^{\{1\}}(y) = D(y)$ and $G^{\{2\}}(y) = G(y) - D(y)$
where $D(y)$ is the newly added operator; unlike Rosenbrock methods, this operator can be nonlinear.  The EIN approach can be used to add and subtract nonlinear operators with time-lagged pieces, however, for reasons of brevity we do not explore this type of EIN method in this work.

%% file: section-motivation.tex
\section{Motivating nonlinear partitioning}\label{sec:motivation}

We now briefly motivate the nonlinearly partitioned ODE \cref{eq:np-ivp} using the viscous Burgers' equation. In particular, we propose three common partitions of the discretized equation, and then show how an integrator that exploits a nonlinear partitioning achieves excellent accuracy and stability without requiring nonlinear solves or the right-hand-side Jacobian. More advanced examples of nonlinear partitioning can be found in \cite{Sebastiano.2023,Boscarino.2016}.

As we have not yet introduced NPRK methods, we only consider the IMEX-NPRK Euler method
\begin{align}
	y_{n+1} = y_n + hF(y_{n+1}, y_n)	.
	\label{eq:nprk-euler}
\end{align}
This integrator treats the first argument of $F(u,v)$ implicitly and the second argument explicitly, and is a nonlinear generalization of the classical IMEX Euler method for additive equations. 
Depending on the choice of $F(u,v)$, the implicit solve at each timestep can be linear or nonlinear. Morover, it follows from Taylor expansion that \cref{eq:nprk-euler} has local truncation error of $\mathcal{O}(h^2)$, analogous to other Euler methods. We also note that \cref{eq:nprk-euler} has appeared previously in the literature, both in the context of time-integration (e.g. the SIRK method from \cite{Boscarino.2016} using an IMEX Euler tableau) and as a first-order method that time-lags certain terms in a nonlinear equation (e.g. this is the standard approach in TRT dating back many decades \cite{Larsen.1988}).

Now, consider the viscous Burger's equation
\begin{align}
    u_t = \epsilon u_{xx} + uu_x
    \label{eq:burgers}
\end{align}
equipped with homogeneous Dirichlet boundary conditions. After spatial discretization with classical second-order finite differences, the resulting ODE system is
\begin{equation}\label{eq:discrete-burgers-1}
    y' = \epsilon D y + \textnormal{diag}(y)Ay,
\end{equation}
for discrete diffusion and advection matrices $D$ and $A$, respectively. A classical explicit or implicit method will treat the equation right-hand-side monolithically. The equation also naturally admits an additive splitting that separates diffusion from the nonlinear term. There are many possibilities for nonlinear splittings; here we  consider one such splitting that is linear in each variable. These three splittings are:
\begin{align}
	\text{Monolithic \cref{eq:generic-ivp}:} &&&  G(y) = \epsilon D y + \textnormal{diag}(y)Ay,\label{eq:burgers-splitting-classical} \\
	\text{Additive \cref{eq:additive-ivp}:} &&& 
		G^{\{1\}}(y) = \epsilon D y, \quad
		G^{\{2\}}(y) = \textnormal{diag}(y)Ay, 
		\label{eq:burgers-splitting-additive} \\
	\text{Nonlinear \cref{eq:np-ivp}:} &&& 
		F(u, v) = \epsilon D u + \text{diag}(v) A u.
		\label{eq:burgers-splitting-nonlinear}
\end{align}
We can also rewrite \cref{eq:discrete-burgers-1} as a second additive system using $G^{\{1\}}(y) = J_n y$, $G^{\{2\}}(y) = G(y) - J_n y$ and the Jacobian $J_n = \frac{\partial G}{\partial y}(y_n) = \epsilon D + \text{diag}(Ay_n) + \text{diag}(y_n)A$; this will lead to a Rosenbock method.
The first-order integrators associated with each of these four splittings are:
\begin{center}
	\vspace{0.5em}
	\renewcommand*{\arraystretch}{1.25}
	\begin{tabular}{l|ll}
		\hline
		Family & ODE & Euler Method 	\\ \hline
		Explicit	& $\dot{y} = G(y)$ & $y_{n+1} = y_n + h G(y_{n})$	\\
		Implicit	& $\dot{y} = G(y)$ & $y_{n+1} = y_n + h G(y_{n+1})$ 	\\
		Rosenbrock & $\dot{y} = G(y)$ & $(I - h J_n)(y_{n+1}-y_n) = h G(y_n)$\\
		
		IMEX		& $\dot{y} = G^{\{1\}}(y) + G^{\{2\}}(y)$ & $y_{n+1} = y_{n} + h G^{\{1\}}(y_{n+1}) + hG^{\{2\}}(y_n)$	\\
		NPRK		& $\dot{y} = F(y, y)$ & $y_{n+1} = y_{n} + h F(y_{n+1},y_n)$	\\ \hline		
	\end{tabular}
	\vspace{0.5em}
\end{center}

\begin{figure}
	
	\centering
	\includegraphics[width=.9\textwidth,trim={0.5cm 4.1cm 1cm 0},clip]{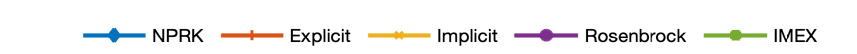}
	
	\includegraphics[width=0.48\textwidth]{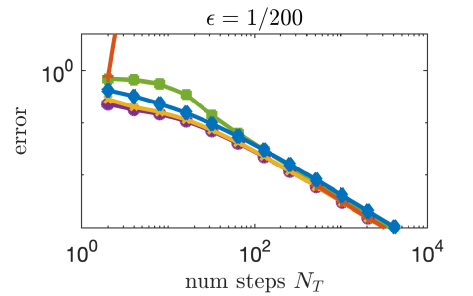}
	\hfill
	\includegraphics[width=0.48\textwidth]{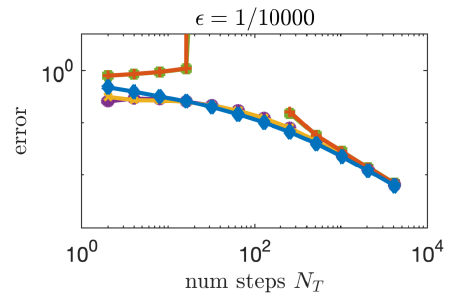}

	\caption{Convergence diagrams for the viscous Burgers' equation \cref{eq:burgers}. We discretized the domain $x \in [-2,2]$ using 1000 spatial grid points, and integrated to time $t=6/10$ using the initial condition $u(x,t=0) = e^{-3x^2}$.}
	\label{fig:burgers-experiment}	
\end{figure}

Next, we compare these five integrators on \cref{eq:discrete-burgers-1} with $\epsilon = 1/200$ and $\epsilon = 1/10000$; \cref{fig:burgers-experiment} contains convergence diagrams. We briefly remark on each method:
\begin{itemize}[leftmargin=2em]
	\item {\em Explicit}. No solves are required at each timestep, however, the largest stable $h$ is limited. For $\epsilon \lessapprox h_x/|u|$ the stepsize is limited by the standard explicit hyperbolic CFL, while for large $\epsilon$ values or small grid spacing $h_x$, diffusion requires that $h = \mathcal{O}(\epsilon/h_x^2)$.
	\item {\em Implicit}. A nonlinear solve is required at each timestep, however the method retains stability at large timesteps regardless of $\epsilon$.
	\item {\em Rosenbrock}. Stability and accuracy are comparable to the implicit method, despite only requiring a linear solve involving the exact right-hand-side Jacobian.
	\item {\em IMEX}. For $\epsilon \ne 0$, a linear implicit solve is required at each step, however the method reduces to explicit Euler when $\epsilon \to 0$. Therefore, improvements in stability compared to explicit Euler require large epsilon or small grid spacing.
	\item {\em NPRK}. The method remains linearly implicit for all $\epsilon$, and provides similar accuracy and stability as implicit and Rosenbrock Euler, while only requiring a linear solve at each step and avoiding the exact right-hand-side Jacobian.
\end{itemize}

We have now seen how nonlinear partitioning can be beneficial for solving a simple one-dimensional PDE. In the next section we introduce a framework that naturally generalizes the method \cref{eq:nprk-euler} by allowing for additional stages, higher-order accuracy, and different types of implicitness.

%% file: section-nprk.tex
\section{Nonlinearly partitioned Runge-Kutta methods}\label{sec:gsi}

We propose a new family of nonlinearly partitioned Runge-Kutta (NPRK) methods for solving the initial value problem
\begin{align}
	y' = F(y,y), \quad y(t_0) = y_0. 
	\label{eq:ivp-nonlinearly-partitioned}	
\end{align}
The ansatz for an $s$-stage NPRK method is
\begin{align}
	\begin{aligned}
		Y_i &= y_n + h\sum_{j=1}^{s} \sum_{k=1}^{s} a_{ijk} F(Y_j,Y_k), \quad i=1,\ldots, s, \\
		y_{n+1} &= y_n + h\sum_{i=1}^s \sum_{j=1}^s b_{ij} F(Y_i,Y_j),	
	\end{aligned}
	\label{eq:nprk-general}
\end{align}
where $Y_i$ are stage values, the rank three tensor $a_{ijk}$ takes the place of the classical Runge-Kutta matrix $a_{ij}$, and the matrix $b_{ij}$ replaces the classical weight vector $b_i$. This additional generality implies that each stage value may be expressed in terms of all possible combinations of stage values substituted into the arguments of $F(u,v)$.

\subsection{Connections to classical RK and ARK} We briefly describe two connections between NPRK methods and existing RK methods. Firstly, the NPRK method \cref{eq:nprk-general} is equivalent to a classical RK method if $a_{ijk}$ is diagonal in the last two indices ($a_{ijk} = 0$ if $j\ne k$), and $b_{ij}$ is diagonal ($b_{ij} = 0$ if $i \ne j$). For more general coefficients, we can interpret an NPRK method as two nonlinearly-coupled classical RK methods. We formalize this concept with the following definition.

\begin{definition}[Underlying method]
	A classical RK integrator $(a_{ij}, b_i, c_i)$ is an underlying method of an NPRK integrator, if \cref{eq:nprk-general} reduces to $(a_{ij}, b_i, c_i)$ when the function $F(u,v)$ in \cref{eq:ivp-nonlinearly-partitioned} depends only on a single argument.	
\end{definition}
The NPRK method \cref{eq:nprk-general} has two underlying RK integrators. If ${F(u,v) = G(u)}$, then it reduces to a classical RK method with coefficients
\begin{align}
   		a^{\{1\}}_{ij} = \sum_{k=1}^s a_{ijk}, \quad  
   		b^{\{1\}}_{i} = \sum_{j=1}^s b_{ij}, \quad
   		c^{\{1\}}_i = \sum_{j=1}^{s} a^{\{1\}}_{ij}.
    \label{eq:nprk-underlying-rk1}
\end{align}
Similarly, if ${F(u,v) = G(v)}$, then \cref{eq:nprk-general}  reduces to a classical RK method with
\begin{align}
       	a^{\{2\}}_{ik} = \sum_{j=1}^s a_{ijk}, \quad
    	b^{\{2\}}_{j} = \sum_{i=1}^s b_{ij}, \quad
    	c^{\{2\}}_i = \sum_{j=1}^{s} a^{\{2\}}_{ij} = c_i^{\{1\}}.
    \label{eq:nprk-underlying-rk2}
\end{align}

A connection also exists between NPRK methods and two-component ARK methods \cite{Kennedy.2003tv4} with shared abscissa. Specifically, if \cref{eq:nprk-general} is applied to the additively partitioned right-hand-side 
	\begin{align}
	y' = F(y,y), \quad 	 F(u,v)=G^{\{1\}}(u) +G^{\{2\}}(v).
	\label{eq:nprk-underlying-ark-equation}
\end{align} 
then it reduces to the ARK method  \cref{eq:nprk-underlying-ark} whose coefficients are those of the underlying methods \cref{eq:nprk-underlying-rk1,eq:nprk-underlying-rk2}. Therefore, NPRK methods generalize ARK methods with shared abscissa vectors in the sense that they solve nonlinearly partitioned problems, but reduce to their underlying ARK method \cref{eq:nprk-underlying-ark,eq:nprk-underlying-rk1,eq:nprk-underlying-rk2} on additively partitioned right-hand-sides. These features have implications for order conditions. Specifically, any NPRK method that is order $q$ accurate, must have underlying methods that are also $q$ accurate and whose coefficients satisfy ARK order conditions up to order $q$. We note that the shared abscissa do not reduce the total number of order conditions for an NPRK methods since order must be mantained on nonlinear partitions (and not exclusively additive ones).

	\subsection{Connections to SIRK}
	
	Any SIRK method (see \cref{subsec:sirk}) can be expressed as an NPRK method with alternating explicit and implicit stages. Specifically, an $s$-stage SIRK method with implicit tableau (${A^I, b}$) and explicit tableau (${A^E, b}$) can be expressed as a $2s$-stage NPRK method with coefficients:
		\begin{align}
				\begin{aligned}
					a_{2i-1,2j,2j-1}	 &= A^E_{ij} \\
					a_{2i,2j,2j-1} &= A^{I}_{ij}	
				\end{aligned} &&
			b_{2i,2i-1} &= b_{i}
		\end{align}
	for $i,j = 1, \ldots s$, and all other $a_{ijk}$ zero. In contrast to SIRK, the NPRK ansatz naturally accommodates more varied types of implicitness, and provides a simple way to derive high-order methods with fewer stages since there is no need to alternate between explicit and implicit stages.

\subsection{Degrees of implicitness}

The sparsity of the NPRK tensor $a_{ijk}$ determines the type of implicitness found in the method. The ansatz \cref{eq:nprk-general} allows for different types of implicitness for each of the two components of the right-hand-side. This includes methods with {\em coupled implicitness} that have implicit stages in both arguments of $F(u,v)$ and methods with {\em decoupled implicitness} for which the unknown for an implicit stage is only found in a single argument of $F(u,v)$. The properties for three types of NPRK methods are described in \cref{tab:nprk-implicitness-sparsity}.
\begin{table}[h]
	\caption{Three example families of NPRK methods.}
	\label{tab:nprk-implicitness-sparsity}
	
	\centering
	\renewcommand*{\arraystretch}{1.5}
	\begin{tabular}{lll}
		First Component & Second Component & Required Sparsity \\ \hline	
		diagonally-implicit  & diagonally-implicit  & $a_{ijk} = 0$ for $j > i$, $k > i$ \\
		explicit & explicit & $a_{ijk} = 0$ for $j\ge i$, $k\ge i$ \\
		diagonally-implicit  & explicit & $a_{ijk} = 0$ for $j > i$, $k \ge i$
	\end{tabular}	
\end{table}

In this paper, we will focus on methods with decoupled implicitness where the first component is treated diagonally-implicitly and the second component is treated explicitly (third option in \cref{tab:nprk-implicitness-sparsity}). Such methods are the nonlinear partitioned counterparts to classical IMEX methods that assume an additively partitioned right-hand-side. We will also derive one class of decoupled implicit-implicit (IMIM) methods that treat both components diagonally implicitly. To reduce the number of free parameters and simplify the structure of the diagonally-implicit coupling we further assume that $a_{i,i,k} = 0$ for all $k \ne i-1$.
The resulting ansatz is
\begin{align}
	\begin{split}
		Y_i &= y_n + h \bigg[ a_{i,i,i-1} F(Y_i,Y_{i-1}) + a_{i,i-1,i} F(Y_{i-1},Y_{i}) + \\
            & \hspace{8ex}\sum_{j=1}^{i-1} \sum_{k=1}^{i-1} a_{ijk} F(Y_j,Y_k) \bigg], \quad i=1,\ldots, s, \\
		y_{n+1} &= y_n + h \sum_{j=1}^s \sum_{k=1}^s b_{jk} F(Y_j,Y_k),
	\end{split}
	\label{eq:nprk-general-imim-imex}
\end{align}
where either $a_{i,i,i-1} = 0$ or $a_{i,i-1,i} = 0$
 depending on whether the stage $Y_i$ is implicit in the first or second argument of $F(u,v)$, respectively. Also note that the constants $a_{1,1,0}$ and $a_{1,0,1}$ are not formally part of the tensor $a_{ijk}$ and are always taken to be zero. 

\subsection{Equivalence with PRK and the resulting order conditions}\label{sec:gsi:order}

Order conditions for a general NPRK method \cref{eq:nprk-general} can be directly obtained by performing rooted tree analysis on the equation $y'=F(u,v)$, and is the focus of a companion paper \cite{nprk2}. In this paper, we will instead show that the NPRK subclass \cref{eq:nprk-general-imim-imex} can be expressed as a PRK integrator \cref{eq:partitioned-runge-kutta}, and then exploit this relationship to obtain order conditions. This approach simultaneously reveals connections with well-known integrators and allows us to reuse existing results \cite{Jackiewicz.2000,Hairer.1993,Kennedy.2003tv4} in the nonlinearly partitioned context.

To show the connection between \cref{eq:partitioned-runge-kutta} and \cref{eq:nprk-general-imim-imex}, we first introduce a subclass of \cref{eq:nprk-general-imim-imex}, named sequentially-coupled NPRK methods. We then show that any method of the form \cref{eq:nprk-general-imim-imex} can be expressed as a sequentially-coupled method, albeit with an increased number of stages. Lastly, it follows that all sequentially-coupled methods can be written as PRK methods \cref{eq:partitioned-runge-kutta}. Moreover, the coefficients of the PRK method are those of the reduced underlying methods of the sequentially-coupled NPRK method.

\begin{definition}
A sequentially coupled method only includes right-hand-side evaluations whose first and second arguments are sequential stages (i.e. $F(Y_j,Y_{j-1})$) and is characterized by $a_{ijk} = 0$ if $k \ne j - 1$, and $b_{ij} = 0$ if $j \ne i-1$.
\end{definition}
The ansatz for a diagonally-implicit, sequentially-coupled, NPRK method is
\begin{align}
	\begin{split}
		Y_1 &= y_n, \\
		Y_i &= y_n + h \bigg[ a_{i,i,i-1} F(Y_{i},Y_{i-1}) + a_{i,i-1,i} F(Y_{i-1},Y_{i}) +
            \\&\hspace{8ex}\sum_{j=2}^{i-1} a_{i,j,j-1} F(Y_{j},Y_{j-1}) \bigg], \quad i = 2, \ldots, s, \\
		y_{n+1} &= y_n + h \sum_{i=2}^s b_{i,i-1} F(Y_{i}, Y_{i-1}).
	\end{split}
	\label{eq:nprk-simple-coupling-imim-imex}
\end{align}
where either $a_{i,i,i-1} = 0$ or $a_{i,i-1,i} =0$. 
We now demonstrate the relationships
\begin{align*}
	\text{Method \cref{eq:nprk-general-imim-imex}} 
	\quad \xleftrightarrow{} \quad
	\text{Method \cref{eq:nprk-simple-coupling-imim-imex}}
	\quad \xrightarrow{} \quad
	\text{PRK \cref{eq:partitioned-runge-kutta}}.
\end{align*}
Since the order conditions for the PRK method \cref{eq:partitioned-runge-kutta} are known \cite{Jackiewicz.2000,Hairer.1993,Kennedy.2003tv4} the equivalence provides order conditions for \cref{eq:nprk-general-imim-imex}.

\begin{proposition}
	The $s$-stage method \cref{eq:nprk-general-imim-imex} is equivalent to an $\hat{s}$-stage sequentially-coupled method \cref{eq:nprk-simple-coupling-imim-imex} with $\hat{s} \ge s$. 
\end{proposition}
\begin{proof}
	See \cref{app:fully-coupled-to-simply-coupled-proof}; the proof also describes a procedure for assembling the coefficient tensor.
\end{proof}

\begin{proposition}
	The $s$-stage stage sequentially-coupled method \cref{eq:nprk-simple-coupling-imim-imex} can be expressed as an $(s-1)$-stage PRK method \cref{eq:partitioned-runge-kutta} with coefficients
\begin{align}
		a_{ij} &= a_{i+1,j+1,j} & 
		b_i&= b_{i+1,i} & 
		c_i&=\sum_{j=1}^{s-1} a_{ij} 
		\label{eq:prk-implicit-method} \\ 
		\widehat{a}_{ij} &= a_{i,j+1,j} &
		\widehat{b}_j &= b_{i+1,i} & 
		\widehat{c}_i &=\sum_{j=1}^{s-1} \widehat{a}_{ij}.	
		\label{eq:prk-explicit-method}
\end{align}
applied to the partitioned equation $z' = F(z,w)$, $w' = F(z,w)$ with $z_n=w_n=y_n$. Note that the weight vectors are identical (i.e. $b = \hat{b}$), and the Tableau of the two RK methods \cref{eq:prk-implicit-method,eq:prk-explicit-method} are related by a simple shift, as visualized in \cref{tab:prk-implicit-explicit-tabalue} for $s=4$.
\end{proposition}
\begin{proof}
	Substitute the coefficients \cref{eq:prk-implicit-method,eq:prk-explicit-method} into \cref{eq:partitioned-runge-kutta}. It follows immediately that the stages and outputs of the PRK method \cref{eq:partitioned-runge-kutta} and NPRK method \cref{eq:nprk-simple-coupling-imim-imex} are related according to $Z_i = Y_{i-1}$, $W_i = Y_{i}$, and $z_{n+1} = w_{n+1} = y_{n+1}$.	
\end{proof}

\begin{table}[h]
	\caption{Tablaue for RK methods \cref{eq:prk-implicit-method,eq:prk-explicit-method} for $s=4$}
	\label{tab:prk-implicit-explicit-tabalue}
	
	\begin{align}
		& (a,b,c): \quad
		\begin{tabular}{r|lll}
			& 0 \\
			$(a_{221})$ & $a_{221}$ \\ 
			$(a_{221} + a_{332})$ & $a_{321}$ & $a_{332}$	\\ \hline
			& $b_{21}$ & $b_{32}$ & $b_{43}$
		\end{tabular}
		\\[1em]
		& (\widehat{a},\widehat{b},\widehat{c}): \quad 
		\begin{tabular}{r|lll}
			$(a_{221})$ & $a_{221}$ \\
			$(a_{221} + a_{332})$ & $a_{321}$ & $a_{332}$	\\
			$(a_{421} + a_{432} + a_{443})$ & $a_{421}$ & $a_{432}$ & $a_{443}$	\\ \hline
			& $b_{21}$ & $b_{32}$ & $b_{43}$
		\end{tabular}
	\end{align}
	
\end{table}

The order conditions required to derive PRK methods \cref{eq:partitioned-runge-kutta} up to order three are contained in \cref{tab:PRK-order-conditions}; they assume that $b = \hat{b}$ which is always true when starting from \cref{eq:nprk-simple-coupling-imim-imex}. For second-order NPRK methods \cref{eq:nprk-simple-coupling-imim-imex}, we only require that the two methods ($a_{ij}$, $b_{i}$, $c_{i}$) and ($\widehat{a}_{ij}$, $\widehat{b}_{i}$, $\widehat{c}_{i}$) are each second order accurate, while for third-order, there are three additional coupling conditions.  See \cite{Jackiewicz.2000,Hairer.1993,Kennedy.2003tv4} for higher-order conditions.

\begin{remark}
\label{remark:oc-nprk-underlying}
 The order conditions of a sequentially-coupled NPRK \cref{eq:nprk-simple-coupling-imim-imex} can also be expressed in terms its underlying methods \cref{eq:nprk-underlying-rk1,eq:nprk-underlying-rk2}. Specifically, the reduced underlying methods (i.e. unused stages removed) are equivalent to the methods \cref{eq:prk-implicit-method,eq:prk-explicit-method} (see \cref{app:sequentially-coupled-underlying-method} for additional details).  
\end{remark}

\begin{table}[h]
	\caption{Order conditions (up to order three) for a PRK method \cref{eq:partitioned-runge-kutta} with shared weight vectors $b_i = \widehat{b}_i$.}
	\label{tab:PRK-order-conditions}
	\begin{align}
		\text{\bf order 1:\quad}& \sum_{i=1}^{s} b_i = 1 \\
		\text{\bf order 2:\quad}& \sum_{i=1}^{s} b_i c_{i} = \sum_{i=1}^{s} b_i \widehat{c}_{i} = \tfrac{1}{2} \\
		\text{\bf order 3:\quad}& \left\{ 
		\begin{aligned}
			& 
				\sum_{i=1}^{s} b_i \left( c_{i} \right)^2 = 
				\sum_{i=1}^{s} b_i c_{i} \widehat{c}_{i} = 
				\sum_{i=1}^{s} b_i \left( \widehat{c}_{i} \right)^2 = 
				\tfrac{1}{3} \\
				&
				\sum_{i=1}^{s} \sum_{j=1}^{s} b_i a_{ij} c_{j} = 
				\sum_{i=1}^{s} \sum_{j=1}^{s} b_i a_{ij} \widehat{c}_{j} = \tfrac{1}{6} \\
				& \sum_{i=1}^{s} \sum_{j=1}^{s} b_i \widehat{a}_{ij} c_{j} = 
				\sum_{i=1}^{s} \sum_{j=1}^{s} b_i \widehat{a}_{ij} \widehat{c}_{j} =
				\tfrac{1}{6}	
		\end{aligned} \right.
	\end{align}
\end{table}

\subsection{Linear stability}\label{sec:gsi:stability}

Classical linear stability analysis \cite{Hairer.1993}[IV.2] studies the behavior of a time integrator applied to the Dahlquist equation $y'=\lambda y$. For Runge-Kutta methods, the timestep iteration reduces to the recurrence $y_{n+1}=R(z)y_n$ where $z=h\lambda$ and the function $R(z)$ is known as the stability function. The aim of linear stability analysis is to determine which $z$ values lead to a bounded iteration. There are many desirable properties for methods including A-stability, which implies that $|R(z)| \le 1$ for all $z \in \mathbb{C}^{-}$, and L-stability which adds the additional condition $\lim_{|z|\to \infty}R(z) = 0$.

In this section, we discuss the linear stability properties of an NPRK method \cref{eq:nprk-general-imim-imex} using the partitioned Dahlquist equation
	\begin{align}
		y' = f(y,y) 
		\quad \text{for} \quad
		f(u,v) = \lambda_1 u + \lambda_2 v.
		\label{eq:dahlquist-partitioned}	
	\end{align}
If rewritten as an additive system \cref{eq:additive}, the scalar equation \cref{eq:dahlquist-partitioned} has been used to study the stability of additive integrators \cite{Ascher.1995,Izzo.2017,Kennedy.2003tv4,Sandu.2015}; in the context of this paper it represents a simplification of a nonlinearly partitioned system. In particular, any system ${y'=F(y,y)}$ can be locally approximated by decoupled partitioned Dahlquist equations if we linearize and assume the resulting component Jacobians are simultaneously diagonalizable (see \cref{app:np-system-dalquist-relation}). 

Since $f(y,y)$ in \cref{eq:dahlquist-partitioned} is additively partitioned, an NPRK method \cref{eq:nprk-general} applied to this problem will reduce to its underlying ARK method \cref{eq:nprk-underlying-ark,eq:nprk-underlying-rk1,eq:nprk-underlying-rk2}. Moreover, the timestep computation reduces to the iteration
\begin{align}
	y_{n+1} = R(z_1,z_2), \quad z_1 = h z_1, \quad z_2 = h z_2,	
\end{align}
with stability function
\begin{align}
	\begin{aligned}
		R(z_1,z_2) &= \frac{\det(\mathbf{I} - (z_1 \mathbf{A}^{\{1\}} + z_2 \mathbf{A}^{\{2\}}) + \mathbf{e}(z_1 \mathbf{b}^{\{1\}} + z_2 \mathbf{b}^{\{2\}})^\text{T} )}{\det(\mathbf{I} - z_1 \mathbf{A}^{\{1\}} - z_2 \mathbf{A}^{\{2\}})}.
	\end{aligned}
	\label{eq:ark-stability-function}
\end{align}
The function $R(z_1,z_2)$ is the stability function for an ARK method \cite{Kennedy.2003tv4} with coefficients 
$\mathbf{A}^{\{1\}}=[a^{\{1\}}_{ij}]$, $\mathbf{b}^{\{1\}}=[b^{\{1\}}_{i}]$, $\mathbf{A}^{\{2\}}=[a^{\{2\}}_{ij}]$, and $\mathbf{b}^{\{1\}}=[b^{\{2\}}_{i}]$ from \cref{eq:nprk-underlying-rk1,eq:nprk-underlying-rk2}. The stability region of the method contains the $(z_1,z_2)$ pairs that lead to a bounded iteration and is defined by
\begin{align}
	\mathcal{S} = \left\{ (z_1,z_2) \in \mathbb{C} \times \mathbb{C} : \left| R(z_1,z_2) \right| \le 1 \right\}.
	\label{eq:ark-stability-region}	
\end{align}
Generally speaking, methods with large stability regions are preferable. However, the types of attainable stability characteristics depend on a method's implicitness, with more implicitness allowing for more stringent stability properties. In the following subsections, we describe desirable stability properties for IMIM and IMEX NPRK methods.

\subsubsection{Stability characteristics for IMIM}

For NPRK methods that are implicit in both arguments of $F(u,v)$, it is reasonable to extend the classical notion of A and L stability. Specifically, a method is A-stable or L-stable if:
\begin{align}
	\label{eq:imim-a-stability}	
	\text{A-stable} & \iff \mathcal{S} \supseteq \mathbb{C}^- \times \mathbb{C}^- \\
	\label{eq:imim-l-stability}	
	\text{L-stable} & \iff \text{A-stable and $\textstyle \lim_{|z_1| \to \infty} R(z_1,z_2) = \lim_{|z_2| \to \infty} R(z_1,z_2) = 0$}.
\end{align}
Similar to classical RK methods, an A-stable NPRK method will retain stability for any pair of $(z_1,z_2)$ with negative real-part, and an L-stable method also achieves maximal damping as either or both of the variables grow arbitrarily large in magnitude.

\subsubsection{Stability characteristics for IMEX}

NPRK methods that are only implicit in one variable cannot be A-stable since either $R(0,z_2)$ or $R(z_1,0)$ will be the stability function of an explicit RK method. However, nothing prevents a method from achieving stability as the implicit variable becomes increasingly stiff. Specifically, we say that a method is A-stable or L-stable in the stiff $z_j$ limit, for $j \in \{1,2\}$, if:
\begin{align}	
	\label{eq:imex-a-stability-zj-limit}
	\text{A-stable in stiff $z_j$ limit} \iff \lim_{|z_j|\to \infty} R(z_1,z_2) \le 1, \\
	\label{eq:imex-l-stability-zj-limit}
	\text{L-stable in stiff $z_j$ limit} \iff \lim_{|z_j|\to \infty} R(z_1,z_2) = 0.
\end{align}

We are also interested in methods that retain stability in the limit as both variables become increasingly stiff. If we let $z_2$ represent the non-stiff variable, then we can assume that $|z_2|<|z_1|$. It is convenient to parametrize $z_2$ in terms of $z_1$, namely $z_2 = \epsilon z_1$ for $|\epsilon|<1$, and consider the reparametrized stability function 
\begin{align}	
	\beta(z_1,\epsilon) = R(z_1,\epsilon z_1).
\end{align}
We then consider the stiff limit as $z_1$ tends to infinity
\begin{align}
	\beta_{\infty}(\epsilon) = \lim_{|z_1|\to \infty} |\beta(z_1,\epsilon)|,
	\label{eq:stiff-z2-beta-inf-def}
\end{align}
and say that a method is stable in the coupled stiff $z_2$ limit if
\begin{align}	
	\beta_{\infty}(\epsilon) \le 1 \text{ for all } |\epsilon| \le 1.
	\label{eq:stiff-z2-condition}
\end{align}

\begin{remark}
The function $\beta_\infty(\epsilon)$ is a polynomial in $\epsilon$ and like all analytic functions, its maximum in a closed region must occur on the boundary (maximum modulus principle). Therefore, a method is stable in the coupled stiff $z_2$ limit if 
\begin{align}	
	\gamma(\theta) = \left| \beta_\infty(e^{i \theta})\right|^2 
	\quad \text{satisfies} \quad
	\gamma(\theta) \le 1 
	\quad \text{\normalfont{for}} \quad \theta \in [0, 2\pi).
	\label{eq:stiff-z2-condition-simplified}
\end{align}
Moreover, for real-valued NPRK coefficients $\gamma(\theta)$ is a cosine series in $\theta$.
For further details, see \cref{app:stability-function-proofs}.
\end{remark}

\subsubsection{Visualizing stability regions}

The four-dimensional stability region \cref{eq:ark-stability-region} is difficult to visualize. One approach is to overlay two-dimensional slices of $\mathcal{S}$ pertaining to a fixed $z_1$ value,
	\begin{align}
		\mathcal{S}(z_1) = \left\{ z_2 \in \mathbb{C} : \left| R(z_1,z_2) \right| \le 1 \right\}.
	\end{align}
	It is also useful to consider the region
	\begin{align}
		\widetilde{\mathcal{S}}(\theta) = \left\{ z_2 \in \mathbb{C} : \max_{\gamma \in \mathbb{R}^+} \left| R(\gamma e^{i\theta},z_2) \right| \le 1 \right\},
	\end{align}
that contains all $z_2$ values that ensure stability for any $z_1=\gamma e^{i\omega}$ in the wedge $\gamma \ge0$ and $\theta \le \omega \le 2\pi - \theta$ (the latter follows from max modulus principle).

When considering stability for PDE discretizations, the values $z_1$ and $z_2$ represent the eigenvalues of the linearized component Jacobians (see \cref{app:np-system-dalquist-relation}). Since real-valued PDE discretizations have spectrums that are symmetric about the imaginary axis, one desires stability for both $z_1$ and $z_1^*$. Therefore, it is more realistic to consider the two-dimensional slices 
\begin{align}
	\mathcal{P}(z_1) &= \left\{ z_2 \in \mathbb{C} : \max 
	\left( 
		\left| R(z_1, z_2) \right|,
		\left| R(z_1^*, z_2) \right|
	\right) \le 1 \right\},
	\label{eq:stability-region-symmetric} \\
	\widetilde{\mathcal{P}}(\theta) &= \left\{ z_2 \in \mathbb{C} : \max_{\gamma \in \mathbb{R}^+} \max_{\omega \in \left\{ \theta, -\theta \right\}}\left| R(\gamma e^{i\omega},z_2) \right| \le 1 \right\},	
	\label{eq:stability-region-symmetric-max-r-theta}	
\end{align} 
that are symmetric about the real $z_1$ axis. When $\arg(z_1) = \pi/2$ these regions approximate stability for a skew-symmetric advection discretization, while $\arg(z_1)=\pi$ approximately represents a symmetric positive-definite diffusion discretization. An intermediate value, $\arg(z_1) \in (\pi/2, \pi)$, approximates a mix of advection and diffusion.		
			
\subsection{$M$-nonlinearly partitioned Runge-Kutta}

The NPRK framework and results contained in this section easily generalize to the $M$-component nonlinearly partitioned equation $y'=F(y,\ldots, y)$ where $F$ has $M$ arguments. The ansatz for an $s$-stage, NPRK$_M$ method is
\begin{align}
	\begin{aligned}
		Y_i &= y_n + h\sum_{j_1,\ldots, j_M=1}^{s} a_{i, j_1,\ldots, j_M} F(Y_{j_1}, \ldots, Y_{j_m}), \quad i=1,\ldots, s, \\
		y_{n+1} &= y_n + h \sum_{j_1,\ldots, j_M=1}^{s} b_{j_1, \ldots, j_M} F(Y_{j_1}, \ldots, Y_{j_M}),	
	\end{aligned}
	\label{eq:nprk-general-m-partitions}
\end{align}
where the tensors $a$ and $b$ are now of rank $M+1$ and $M$, respectively. The method \cref{eq:nprk-general-m-partitions} nonlinearly combines $M$ underlying RK methods and when applied to an additively partitioned system $F(Y_{1}, \ldots, Y_{M}) = \sum_{k=1}^M F^{\{k\}}(Y_{k})$ it reduces to an ARK$_M$ method. For brevity we postpone additional study of NPRK$_M$ methods to future work; order conditions for \cref{eq:nprk-general-m-partitions} are presented in the companion paper \cite{nprk2}.

%% file: section-example-methods.tex
\section{Example NPRK Methods}\label{sec:example}

Now that we have introduced the NPRK framework, we derive several example IMEX and IMIM NPRK methods of order two and three. This section lists the method coefficients and the important steps used in each derivation.  Detailed derivations can be found in \cref{app:method-derivations-imex,app:method-derivations-imim}. 
\subsection{IMEX NPRK}

We are primarily interested in IMEX NPRK methods that are L-stable in the stiff $z_1$-limit (see \cref{eq:imex-l-stability-zj-limit}) and have optimized stability in the coupled stiff $z_2$-limit (see \cref{eq:stiff-z2-condition}). For simplicity we consider only $s$ stage, sequentially-coupled IMEX methods that can be represented using the following two-dimensional tableau:
\begin{align}
	\label{eq:nprk-imex-simply-coupled-with-tableau}
	\begin{small}
		\begin{aligned}
			Y_i &= y_n + h\sum_{j=2}^i a_{i,j,j-1} F(Y_j,Y_{j-1}), \quad i=1,\ldots, s, \\
			y_{n+1} &= y_n + h \sum_{j=2}^s b_{j,j-1} F(Y_j,Y_{j-1}).
		\end{aligned}
	\end{small} && 
	\begin{small}
		\begin{tabular}{l|llll}
			& 0 \\
			& $a_{221}$ & \\
			& $\vdots$ & $\ddots$ & \\ 
			& $a_{s21}$ & $\hdots$ & $a_{s,s,s-1}$ \\ \hline
			& $b_{21}$  & $\hdots$ & $b_{s,s-1}$
		\end{tabular}
	\end{small}
\end{align}
We consider methods with up to five stages and four implicit solves. Our convention is IMEX-NPRK$o$[$sq$] with $o$ representing order, $s$ representing stages, and $q$ representing the number of implicit solves. For certain methods we append the suffixes ``Si'' and ``Sa'' for singly-implicit and stiffly-accurate methods, respectively. Lastly, although it is possible to derive IMEX-NPRK method with underlying implicit integrators of DIRK or ESDIRK type, we focus our attention on the DIRK case.

\subsubsection{IMEX two-stage methods}

The simplest sequentially coupled IMEX-NPRK method is the two-stage, first-order method (first discussed in \cref{sec:motivation})
\begin{align}
	&\begin{aligned}
		Y_{1}   &= y_n \\
		Y_{2}   &= y_n + h F(Y_2,Y_1) \\	
		y_{n+1} &= Y_{2}	
	\end{aligned}
	&&
	\begin{aligned}
		&\text{\underline{IMEX-NPRK1[21] (Euler)}} \\
		&\begin{tabular}{l|lll}
		& 0 \\
		& 1 & \\ \hline
		& 1 &  &
	\end{tabular}
	\end{aligned}
	\label{eq:imex-nprk1-21}
\end{align}
This method is L-stable in $z_1$ and is stable in the coupled stiff $z_2$ limit; its stability functions are $R(z_1,z_2) = \frac{1+z_2}{1-z_1}$ and $\beta_{\infty}(\epsilon) = \epsilon$. 

\subsubsection{IMEX three-stage methods}

Next we consider three-stage, sequentially coupled IMEX-NPRK methods \cref{eq:nprk-imex-simply-coupled-with-tableau} with $s=3$.
This ansatz admits a second-order method that requires a single implicit solve:
\begin{align}
	\renewcommand\arraystretch{1.25}
	\begin{tabular}{l|lll}
		& 0 \\
		& $\tfrac{1}{2}$ & \\
		& $\tfrac{1}{2}$ & 0 & \\ \hline
		& $0$  & $1$
	\end{tabular}
	&& \text{reduces to} &&
	\begin{aligned}
		& \text{\underline{IMEX-NPRK2[31] (Midpoint)}} \\
		& \begin{aligned}
			Y_1 &= y_n \\
			Y_2 &= y_n + \tfrac{h}{2} F(Y_2,Y_1) \\
			y_{n+1} &= y_n + h F(Y_2,Y_2)
		\end{aligned}
	\end{aligned}
	\label{eq:IMEX-NPRK2-31}
\end{align}
This method is reduced to its simpler form by observing that stages $Y_3$ and $Y_2$ of the simply coupled method are identical, thus $Y_3$ can be removed and $F(Y_3,Y_2)$ can be replaced with $F(Y_2,Y_2)$ in the output computation. The underlying implicit and explicit methods are the well-known implicit and explicit midpoint methods, respectively. This method is not L-stable in $z_1$, however it retains stability so long as $|z_2 + 1| \le 1$, and should only be used if $F(u,v)$ is entirely non-stiff in the second argument. The stability functions are
\begin{align}
	R(z_1,z_2) &= \frac{z_1 (z_2+1)+ 1 + (z_2 + 1)^2}{2-z_1}	&
	\lim_{|z_1| \to \infty} R(z_1,z_2) &= -(z_2 + 1)
\end{align}

\begin{figure}

	\centering
	\begin{small}
		\begin{tabular}{cccc}
			& $\text{arg}(z_1) = \pi$ & $\text{arg}(z_1) = \frac{3\pi}{2}$ & $\text{arg}(z_1) = \frac{\pi}{2}$ \\		
			\rotatebox{90}{\scriptsize \hspace{4em} 2[31] -- \cref{eq:IMEX-NPRK2-31}} &
			\includegraphics[align=b,width=0.28\textwidth]{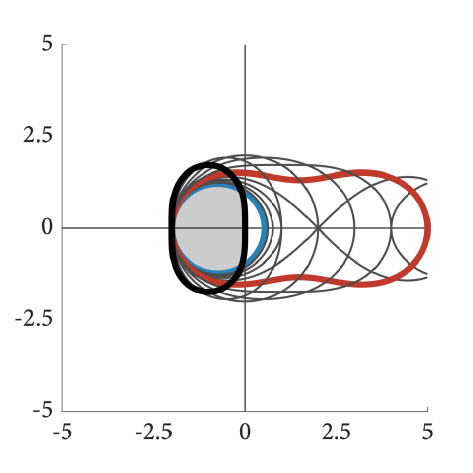} &
			\includegraphics[align=b,width=0.28\textwidth]{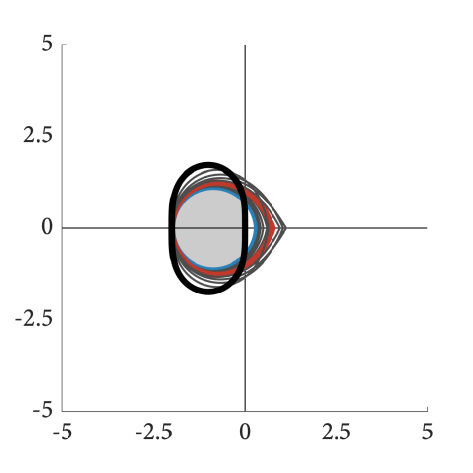} &
			\includegraphics[align=b,width=0.28\textwidth]{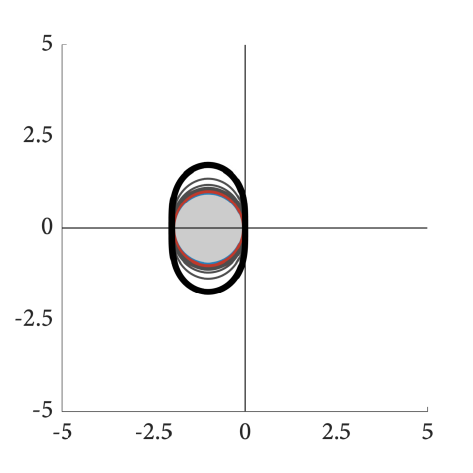} \\
			\rotatebox{90}{\scriptsize \hspace{1em} 2[32]a/2[42]a -- (\ref{eq:NPRK-IMEX-32}/\ref{eq:IMEX-NPRK-42})} &
			\includegraphics[align=b,width=0.28\textwidth]{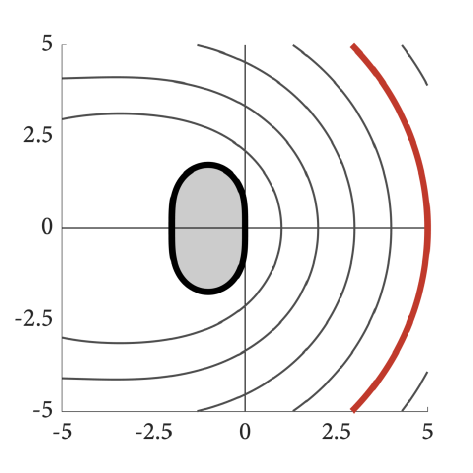} &
			\includegraphics[align=b,width=0.28\textwidth]{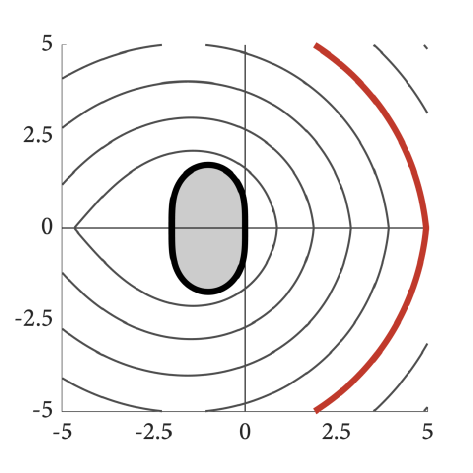} &
			\includegraphics[align=b,width=0.28\textwidth]{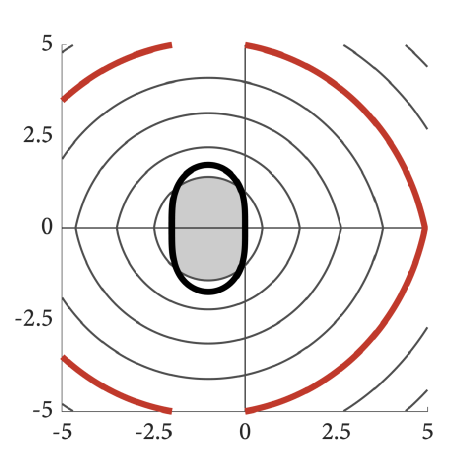} \\
			\rotatebox{90}{\scriptsize \hspace{1em} 2[32]b/2[42]b -- (\ref{eq:NPRK-IMEX-32}/\ref{eq:IMEX-NPRK-42}) } &
			\includegraphics[align=b,width=0.28\textwidth]{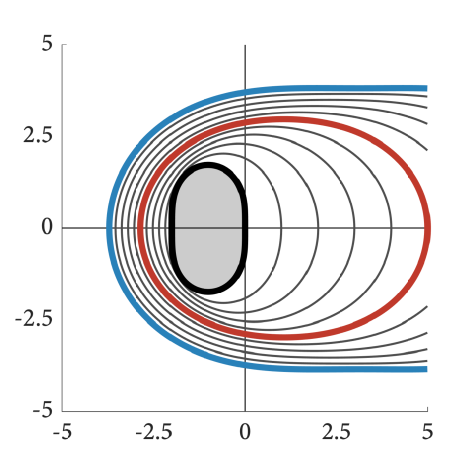} &
			\includegraphics[align=b,width=0.28\textwidth]{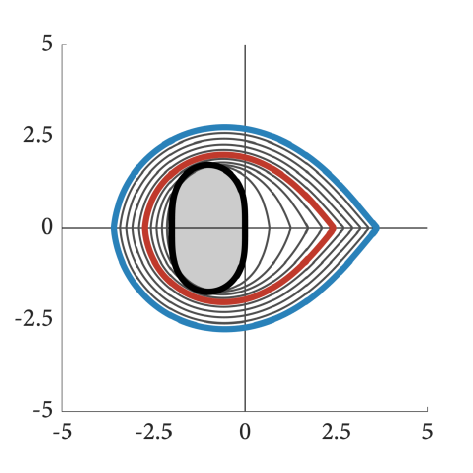} &
			\includegraphics[align=b,width=0.28\textwidth]{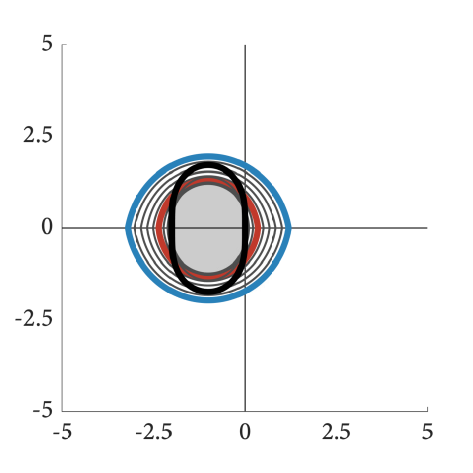}
		\end{tabular}
	\end{small}
		
	\begin{tabular}{c}
		{\tiny \textcolor{black}{\hdashrule[0.2ex]{2em}{2pt}{}} $\left|z_1\right| = 0$ \hspace{1em}}
		{\tiny \textcolor{plot_red}{\hdashrule[0.2ex]{2em}{2pt}{}} $\left|z_1\right| = 5$ \hspace{1em}}
		{\tiny \textcolor{plot_blue}{\hdashrule[0.2ex]{2em}{2pt}{}} $\left|z_1\right| = 10$}
	\end{tabular}

	\caption{Linear stability region slices for IMEX-NPRK 2[31] (top row), 2[32]a or 2[42]a (middle row), and 2[32]b or 2[42]b (bottom row). Each contour is a two-dimensional slice $\mathcal{P}(z_1)$, from \cref{eq:stability-region-symmetric}, with $\arg(z_1) \in \{\pi, \frac{3\pi}{4}, \frac{\pi}{2}\}$ and $|z_1| \in \{0,1,\ldots,10\}$; contours corresponding to $|z_1|=0,5,10$ are colored black, red, and blue, respectively. These choices of $\arg(z_1)$ respectively approximate an implicit linear component with diffusion, a mix of diffusion and oscillation, and pure oscillation. The regions shaded in grey are $\widetilde{\mathcal{P}}(\arg(z_1))$ from \cref{eq:stability-region-symmetric-max-r-theta}.}
	
	\label{fig:second-order-stability-overlays}
	
\end{figure}

Starting from \cref{eq:nprk-imex-simply-coupled-with-tableau} with $s=3$, we can also obtain
a one-parameter family of second-order methods that require two solves and are L-stable in the stiff $z_1$-limit:
\begin{align}
	\begin{aligned}
		a_{221} &= \frac{1}{2 b_{32}}, &
		a_{321} &= \frac{-2 b_{32}^3+6 b_{32}^2-4 b_{32}+1}{2 b_{32}^2 (2 b_{32}-1)}, \\
		a_{332} &= \frac{b_{32}-1}{2 b_{32}-1}, &
		b_{21}  &= 1-b_{32}
	\end{aligned}
\end{align}
The two choices which maximize stability in the coupled stiff $z_2$-limit are:	
	\begin{align}
		\begin{tabular}{ccc}
			\text{\underline{IMEX-NPRK2[32]a}} & \hspace{5em} & \text{\underline{IMEX-NPRK2[32]b}}	\\[0.5em]
			\text{(Global Optimum)} & & \text{(Local Optimum)} \\[0.5em]	
			\renewcommand\arraystretch{1.25}
			\begin{tabular}{l|lll}
				& $0$ \\
				& $(1+\frac{1}{\sqrt{2}})$		\\
				& $(-2-\frac{3}{\sqrt{2}})$	& $(1+\frac{1}{\sqrt{2}})$	 \\ \hline
				& $\frac{1}{\sqrt{2}}$ 		& $(1-\frac{1}{\sqrt{2}})$
			\end{tabular}
			& & \renewcommand\arraystretch{1.275}
			\begin{tabular}{l|lll}
				& $0$ \\
				& $(1-\frac{1}{\sqrt{2}})$		\\
				& $(\frac{3}{\sqrt{2}}-2)$		&	$(1-\frac{1}{\sqrt{2}})$	\\ \hline
				& $-\frac{1}{\sqrt{2}}$ 		& $(1+\frac{1}{\sqrt{2}})$
			\end{tabular}
		\end{tabular}
		\label{eq:NPRK-IMEX-32}
	\end{align}
	NPRK2[32]a is stable in the coupled $z_2$-limit \cref{eq:stiff-z2-condition}, while NPRK2[32]b is not (see \cref{fig:imex-232-stability-extra}). As expected, the stability region of NPRK2[32]a is significantly larger than that of NPRK2[32]b (See \cref{fig:second-order-stability-overlays}). However, NPRK2[32]b has a significantly smaller error constant; compare \cref{eq:nprk232s-a-error,eq:nprk232s-b-error}.

\subsubsection{IMEX four-stage methods} 

Next we consider four-stage, sequentially-coupled IMEX NPRK methods \cref{eq:nprk-imex-simply-coupled-with-tableau} with $s=4$.
This anstaz contains second-order methods with improved error coefficients and up to three implicit solves. We begin by considering methods that only require two solves. A simple family of second-order methods that are L-stable in the stiff $z_1$ limit is
\begin{align}
	\begin{aligned}
		a_{332} &= 0, \quad
		a_{421} = \frac{a_{221} (2 a_{221} (b_{43}-1)-4 b_{43}+3)+b_{43}-1}{2 (a_{221}-1) b_{43}}, \\
		a_{432} &= 0, \quad a_{443} = \frac{1}{2 (a_{221}-1)}+1, \quad
		b_{21} = 1-b_{43}, \quad
		b_{23} = 0.
	\end{aligned}
\end{align}
Optimizing stability in the stiff $z_2$-limit, leads to the choices
\begin{align}
	\text{global optimum: } a_{221} = 1 + \tfrac{1}{\sqrt{2}},	&&
	\text{local optimum: } a_{221} = 1 - \tfrac{1}{\sqrt{2}}.
\end{align}
The remaining free parameter $b_{43}$ is selected to minimize the two-norm of the third-order residual vector. This leads to the methods
\begin{align}
	\begin{tabular}{ccc}
		\text{\underline{IMEX-NPRK2[42]a}} & \hspace{3em} & \text{\underline{IMEX-NPRK2[42]b}}	\\[0.5em]
		\text{(Global Optimum)} & & \text{(Local Optimum)} \\[0.5em]	
		\renewcommand\arraystretch{1.75}
		\begin{tabular}{l|llll}
			& $0$ \\
			& $1+\frac{1}{\sqrt{2}}$		\\
			& $\frac{26 - 3 \sqrt{2}}{42}$	& $0$ \\
			& $\frac{-20 -23 \sqrt{2}}{42}$ & $0$ & $1+\frac{1}{\sqrt{2}}$ \\ \hline 
			& $\frac{16-9 \sqrt{2}}{94}$ 		& $0$ & $\frac{78+9 \sqrt{2}}{94}$ \\ 
		\end{tabular}
		& & \renewcommand\arraystretch{1.75}
		\begin{tabular}{l|llll}
			& $0$ \\
			& $1-\frac{1}{\sqrt{2}}$		\\
			& $\frac{26 + 3 \sqrt{2}}{42}$	& $0$ \\
			& $\frac{-20 + 23 \sqrt{2}}{42}$ & $0$ & $1-\frac{1}{\sqrt{2}}$ \\ \hline 
			& $\frac{16 + 9 \sqrt{2}}{94}$ 		& $0$ & $\frac{78-9 \sqrt{2}}{94}$ \\ 
		\end{tabular}
	\end{tabular}
	\label{eq:IMEX-NPRK-42}
\end{align}
The linear stability regions of the methods \cref{eq:IMEX-NPRK-42} are identical to those of  \cref{eq:NPRK-IMEX-32}. However, these four stage methods have significantly improved error constants at the cost of requiring  the storage of one additional stage; compare \cref{eq:nprk242s-a-error,eq:nprk242s-b-error} to \cref{eq:nprk232s-a-error,eq:nprk232s-b-error}.

Next, we consider four-stage methods with three implicit solves. Imposing second-order conditions leads to the six-parameter method family \cref{eq:four-stage-three-solve-second-order-family}. One possibility for reducing free parameters is to impose the singly-implicit condition $a_{i,i,i-1}= \gamma$ and L-stability in $z_1$. These choices lead to the three parameter family
	\begin{align}
		\begin{small}
		\begin{aligned}
				a_{321} &= \frac{1-2 \gamma  (b_{32}+b_{43})}{2 b_{43}}, \quad
			a_{432} = \frac{\gamma  (-2 (\gamma -2) \gamma -1)}{2 \gamma  (b_{32}+b_{43})-1}, \quad
			b_{21} = -b_{32}-b_{43}+1, \\
			a_{421} &= \frac{1}{2} \left(\frac{b_{32} (2 b_{32} \gamma -1)}{b_{43}^2}+\frac{2 (b_{32}-1) \gamma +1}{b_{43}}+\frac{2 \gamma  (2 (\gamma -2) \gamma +1)}{2 \gamma  (b_{32}+b_{43})-1}\right).
		\end{aligned}
		\end{small}
	\end{align}
	A choice that simultaneously optimizes error and stability in the stiff $z_1$ limit is $b_{32} = -0.0054849$, $b_{43} = 0.237378$ and $a_{221} = 0.553658$. The resulting method is
\begin{align}
	\begin{tabular}{c}
		\text{\underline{IMEX-NPRK2[43]-Si}}	\\[0.5em]	
		\begin{small}
		\renewcommand\arraystretch{1.25}
		\begin{tabular}{l|llll}
			& $0$ \\
			& $0.553658000000000$		\\
			& $1.565480078356882$	& $0.553658000000000$ \\
			& $0.258254432834781$ & $-0.448126696979816$ & $0.553658000000000$ \\ \hline 
			& $0.768106900000000$ & $-0.005484900000000$ & $0.237378000000000$ \\ 
		\end{tabular}
		\end{small}
		\label{eq:IMEX-NPRK2-43-Si}	
	\end{tabular}
\end{align}
Another direction for reducing free parameters in \cref{eq:four-stage-three-solve-second-order-family} is to impose the stiff accuracy condition $b_{i,i-1} = a_{s,i,i-1}$, $i=2,\ldots,4$. This leads to the three parameter family of methods \cref{eq:four-stage-three-solve-second-order-stiff-family}, that is L-stable in $z_1$ and has stability function ${\beta(\epsilon) = -\epsilon^3}$ (this implies the method is stable in the coupled stiff $z_2$ limit). Further imposing a singly-implicit structure  $a_{i,i,i-1}= \gamma$ leads to the one parameter family of methods
\begin{align}
	\begin{small}
		\begin{aligned}
			a_{321} &= \frac{1 -2 \gamma^2 +  f(\gamma)}{4\gamma}, &	
			a_{421} &= \frac{-1 + 4\gamma -2 \gamma^2 + f(\gamma)}{4\gamma},	&
			a_{432} &= \frac{1 - 2 \gamma^2 -  f(\gamma)}{4\gamma} &
		\end{aligned}
	\end{small}
\end{align}
with $f(\gamma) = \sqrt{1-4 \gamma^2 (\gamma (3 \gamma-8)+3)}$. Taking $\gamma = 0.386585$ leads to a method with optimal error; its tableau is
\begin{align}
	\begin{tabular}{c}
		\text{\underline{IMEX-NPRK2[43]-SiSa}}	\\[0.5em]	
		\begin{small}
		\renewcommand\arraystretch{1.25}
		\begin{tabular}{l|llll}
			&			   	   0   &                   &                    \\
			&   0.386585000000000  &                   &                    \\
			&   1.027233588987035  & 0.386585000000000 &                    \\
			&   0.733856970649542  &-0.120441970649542 &  0.386585000000000 \\ \hline 
			&   0.733856970649542  &-0.120441970649542 &  0.386585000000000
		\end{tabular}
		\end{small}
	\end{tabular}
	\label{eq:IMEX-NPRK2-43-SiSa}	
\end{align}

\subsubsection{IMEX five-stage methods}

We require at least five stages ($s=5$) to achieve third-order accuracy using the sequentially-coupled IMEX NPRK \cref{eq:nprk-imex-simply-coupled-with-tableau}. There are many degrees of freedom in such methods and Mathematica is unable to provide close-form solutions. We therefore consider several simplifications. If we restrict ourselves to stiffly accurate methods $b_{i,i+1} = a_{s,i,i+1}$ then there are only two solutions. The first is
\begin{align}
	&\nonumber \text{\underline{IMEX-NPRK3[54]-Sa}} \\
	& \renewcommand\arraystretch{1.25}
	\begin{tabular}{l|llll}
		& $1$ \\
		& $-\tfrac{2}{3}$ & $\frac{2}{3}$ \\
		& $\tfrac{5}{12}$ & $-\tfrac{5}{12}$ & $\tfrac{1}{2}$ \\
		& $-\tfrac{1}{2}$ & $\tfrac{1}{6}$ & $\tfrac{2}{3}$ & $\tfrac{2}{3}$ \\ \hline
		& $-\tfrac{1}{2}$ & $\tfrac{1}{6}$ & $\tfrac{2}{3}$ & $\tfrac{2}{3}$
	\end{tabular}
	\label{eq:IMEX-NPRK3-54-Sa}
\end{align}
and the second is listed in \cref{eq:third-order-stiffly-accurate-b} (however it possesses an undesirable diagonal coefficient of approximately $4.55$). If we restrict ourselves to singly-diagonally implicit methods $a_{i,i,i-1} = \gamma$, $i=2, \ldots, s$, then one possible choice with decent stability in the stiff $z_2$ limit is:
\begin{align}
& \text{\underline{IMEX-NPRK3[54]-Si ($\gamma = 54/100$)} -- coefficients listed to 20 digits of precision} \label{eq:IMEX-NPRK354-Si} \\
& \nonumber \begin{tiny}
\begin{tabular}{c|cccc}
	 & 0.54000000000000000000 & 0 & 0 & 0 \\
	 & 0.10402085874596586377 & 0.54000000000000000000 & 0 & 0 \\
	 & -1.2409681743028102473 & 0.42383482979738431001 & 0.54000000000000000000 & 0 \\
	 & 0.42903447708369521671 & -1.0829950086155536734 & 0.24651165580639138296 & 0.54000000000000000000 \\ \hline
	 & -0.32058288115984556996 & 1.0095140978756513629 & 0.044585281470753018014 & 0.26648350181344118903
\end{tabular}
\end{tiny}%
\end{align}
We note that this method is not stable in the stiff coupled $z_2$ limit. More generally, there is also a class of methods with arbitrary $\gamma$, however we are unable to obtain closed form coefficients using Mathematica.

\subsection{IMIM NPRK} 

An IMIM method may be desirable when: (1) both arguments of $F(y,y)$ are stiff processes, and (2) there exist efficient implicit solvers for each  argument of $F$ (e.g. solve for $Y$ in $Y = b + h F(Y,x)$ and $Y = b + h F(x,Y)$), but no efficient solvers for the full right-hand-side (e.g. solve for $Y$ in $Y = b + h F(Y,Y)$). We briefly describe two second-order, A-stable IMIM methods; detailed derivations can be found in \cref{app:method-derivations-imim}. Our starting ansatz is
\begin{align}
	\begin{aligned}
		Y_1 &= y_n, \\
		Y_2 &= y_n + h a_{221} F(Y_2,Y_1), \\
		Y_3 &= y_n + h a_{321} F(Y_2,Y_1) + h a_{323} F(Y_2,Y_3), \\
		y_{n+1} &= y_n + hb_{21} F(Y_2,Y_1) + hb_{23} F(Y_2,Y_3).
	\end{aligned}
\end{align}
There is a one-parameter family of methods with bounded stability functions as $z_1$ or $z_2$ tend to infinity; the associated coefficients are
\begin{small}
	\begin{align}
		b_{21}  &= 1-b_{23}, &
		a_{221} &= \frac{1}{2}, & 
		a_{321} &= \frac{1-b_{23}}{4 (b_{23}-1) b_{23}+2}, &
		a_{323} &= \frac{b_{23}}{4 (b_{23}-1) b_{23}+2}.
	\end{align}
\end{small}%
If we take $b_{23} \in \{\tfrac{1}{2}, 1\}$ then the method has a separable stability function 
	\begin{align}
		R(z_1,z_2) = f(z_1)g(z_2) \quad \text{with} \quad f(z) = g(z) = \tfrac{z+2}{z-2}.	
	\end{align}
Since both $f(z)$ and $g(z)$ are both A-stable (i.e. $\max(|f(z)|,|g(z)|) \le 1$ for $\text{Re}(z) \le 0$), then the NPRK method is A-stable (in both $z_1$ and $z_2$). Moreover, both choices of $b_{32}$ lead to well-known underlying integrators. For $b_{32} = 1$, $M^{\{1\}}$ and $M^{\{2\}}$ are both the implicit midpoint method, while for $b_{32} = 1/2$, $M^{\{1\}}$ is the implicit midpoint and $M^{\{2\}}$ is the implicit Crank-Nicolson method.
\begin{align}
	\nonumber
	& \text{\underline{IMIM Midpoint $(b_{23} = 1)$}} &
	& \text{\underline{IMIM Midpoint/Crank-Nicolson $(b_{23} = \tfrac{1}{2})$}} \\
	& \begin{aligned}
		Y_1 &= y_n \\
		Y_2 &= y_n + \tfrac{h}{2} F(Y_2,Y_1) \\
		Y_3 &= y_n + \tfrac{h}{2} F(Y_2,Y_3) \\
		y_{n+1} &= y_n + h F(Y_2,Y_3)
	\end{aligned}
	&&
	\begin{aligned}
		Y_1 &= y_n \\
		Y_2 &= y_n + \tfrac{h}{2} F(Y_2,Y_1) \\
		Y_3 &= y_n + \tfrac{h}{2} F(Y_2,Y_1) + \tfrac{h}{2} F(Y_2,Y_3) \\
		y_{n+1} &= y_n + \tfrac{h}{2} F(Y_2,Y_1) + \tfrac{h}{2} F(Y_2,Y_3)
	\end{aligned}
	\label{eq:imim-second-order}
\end{align}

\begin{remark} 
For IMIM methods one can feely exchange the order of the arguments of F without affecting the accuracy or linear stability; this is equivalent to redefining $F(u,v)$ according to $F(u,v) := F(v,u)$. See \cref{eq:imim-midpoint-with-transpose,eq:imim-crank-with-transpose} for ``flipped'' variants of \cref{eq:imim-second-order}.
\end{remark}

%% file: section-algorithm-imex.tex
\section{A practical algorithm for sequentially coupled IMEX-NPRK}
\label{sec:alg}

In practice it can be convenient to let $F(y,y)$ take the form
\begin{align}
    F(u,v) \coloneqq F_E(v) + F_I(u,v),
    \label{eq:practical_nonlinear_partitioning_f}
\end{align}
with $F_E(v)$ denoting a non-stiff component that will be treated explicitly, and $F_I(u,v)$ a stiff component that will be treated implicitly in the first argument $u$.
This form naturally isolates additive or component-wise partitioned terms that should always be treated in a fully explicit manner; for example, in the coupling of radiation and hydrodynamics \cite{southworth2023implicit}, the hydrodynamic equations are treated explicitly and contained in $F_E(v)$. This approach is more natural than indirectly representing explicitness within an implicit solve over all components. Assuming \cref{eq:practical_nonlinear_partitioning_f}, each stage of a sequentially coupled IMEX-NPRK method takes the form
\begin{align}
    r_{i} & = y_n +
        \Delta t \sum_{j=1}^{i-1} a_{i,j,j-1}F(Y_{j}, Y_{j-1}), \\
    Y_{i} - a_{i,j,j-1}\Delta tF_I(Y_{i}, Y_{i-1})& =
        r_{i} + a_{i,j,j-1}\Delta t F_E(Y_{i-1}).
\end{align}

Altogether, a practical and moderately efficient implementation is given in \Cref{alg:nprk}.\footnote{Note, one can save further memory by reusing vectors carefully, but the resulting algorithm is much less clean.} The implementation is done in terms of stage values, as in many multiphysics codes this is simpler and more natural than using the stage derivative formulation. We also avoid computing the action of the operators as much as possible due to both cost of assembly (e.g., assembling high-order finite element matrices can be fairly expensive), and difficulty of applying the implicit operator in certain cases (e.g., transport codes are often set up to \emph{solve} an implicit transport equation, but not necessarily \emph{apply} it as an operator). For the implicit operator in particular, we compute its action indirectly following the implicit solve.

\begin{algorithm}[!htb]
  \caption{Sequentially coupled IMEX-NPRK time step
    \label{alg:nprk}}
  \begin{algorithmic}[1]
    \Let{$Y_{\text{old}}$}{$y_n$}
    \Let{$Y_{\text{new}}$}{$y_n$}
    \Let{$y_{(n+1)}$}{$y_n$}\Comment{Only needed if not stiffly accurate}
    \Let{$r_{i}$}{$\mathbf{0}$ for $i=1,...,s$}
    \Statex

    \For{$j = 1 \textrm{ to } s$}\Comment{Loop over stages $1,...,s$}

        \Let{$Y_{\text{old}}$}{$Y_{\text{new}}$}
        \If{$a_{j,j,j-1} \neq 0$}\Comment{Standard implicit stage}
            \State{$\boldsymbol{\delta} = \Delta tF_E(Y_{\text{old}})$}
            \Comment{Compute explicit part operator}
            \State{$r_{j} \pluseq y_n + a_{j,j,j-1}\boldsymbol{\delta}$}
            \Comment{Form implicit right-hand side in place}
            \State{Solve $Y_{\text{new}} - a_{j,j,j-1}\Delta t F_I(Y_{\text{new}},Y_{\text{old}}) = r_{j}$}
            \Comment{Implicit solve for $Y_{\text{new}}$}
            \State{$\boldsymbol{\delta} \pluseq (Y_{\text{new}} - r_{j})/ a_{j,j,j-1}$}
            \Comment{Define $\boldsymbol{\delta} = \Delta tF(Y_{\text{new}},Y_{\text{old}})$}

        \Else\Comment{Explicit stage in implicit variable}
            \State{$Y_{\text{new}} = Y_{\text{old}}$} 
            \State{$\boldsymbol{\delta} = \Delta tF(Y_{\text{new}},Y_{\text{old}})$}
        \EndIf
        \Statex{}

        \For{$i = (j+1) \textrm{ to } s$}
            \State{$r_{i} \pluseq {a}_{i,j,j-1}\boldsymbol{\delta}$}
            \Comment{Update future implicit right-hand sides}
            \State{$y_{(n+1)} \pluseq {b}_{i,i-1}\boldsymbol{\delta}$}
            \Comment{Update next solution}
        \EndFor
    \EndFor
    \State\Return $y_{n+1}$
    \end{algorithmic}
\end{algorithm}

%% file: section-numerical-experiments.tex
\section{Numerical experiments}
\label{sec:results}

\subsection{Revisiting Burger's equation}
\label{sec:results:burger}

When nonlinearly partitioning an equation, it is not always straightforward (or possible) to separate stiff and non-stiff terms. For this reason, we desire IMEX-NPRK methods that remain stable even when explicitly treated terms become stiff; this is precisely what stability in the coupled stiff $z_2$ limit \cref{eq:stiff-z2-condition} attempts to achieve. Here, we revisit Burgers' equation \cref{eq:burgers}, to demonstrate the utility of methods with this property.

In \cref{sec:motivation} we neglected the fact that the nonlinearity in Burgers can be written in either conservative or non-conservative form. This leads to two distict discretizations
\begin{align}
	\text{non-conservative:} & \quad u u_x \to \mathbf{u} .* A \mathbf{u} &
	\text{conservative:} & 	\quad \tfrac{1}{2} (u^2)_x \to \tfrac{1}{2} A ( \mathbf{u}.* \mathbf{u})
\end{align}
where $A$ is the discrete advection matrix. For the non-conservative representation, it is straightforward to identify that the Jacobian of the term $A\mathbf{u}$ will have larger eigenvalues than that for the lone $\mathbf{u}$. However, in the conservative representation the two $\mathbf{u}$ are interchangeable. The associated nonlinear partitions for the discretized viscous Burgers' equation are
\begin{align}
	\label{eq:burgers-nonlinear-partitioning-non-conservative}
	\text{non-conservative:} & \quad F(\mathbf{u},\mathbf{v}) = D\mathbf{u} + \text{diag}(\mathbf{v})  A \mathbf{u} \\
	\label{eq:burgers-nonlinear-partitioning-conservative}
	\text{conservative:} & \quad F(\mathbf{u},\mathbf{v}) = D\mathbf{u} + \tfrac{1}{2} A ( \text{diag}(\mathbf{v}) \mathbf{u} ),	
\end{align}
where only the first splitting clearly separates stiff and non-stiff terms (note, the non-conservative $F(\mathbf{u},\mathbf{u})$ is identical to the one considered in \cref{sec:motivation}). Nevertheless, it is typically desirable to treat the nonlinearity in a conservative way. 

\begin{figure}
	
	\centering
	\hspace{0.75em} \includegraphics[width=.95\textwidth,trim={0.5cm 4.1cm 0.4cm 0},clip]{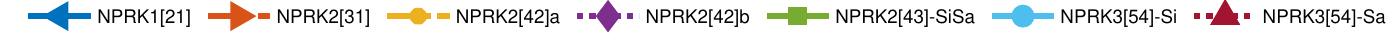}
	
	\vspace{0.5em}
	\includegraphics[width=0.48\textwidth]{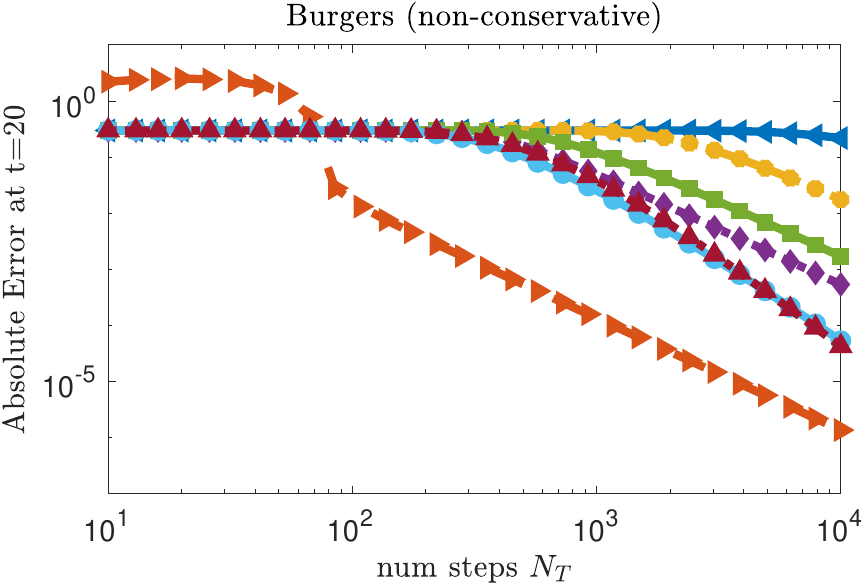}
	\hfill
	\includegraphics[width=0.48\textwidth]{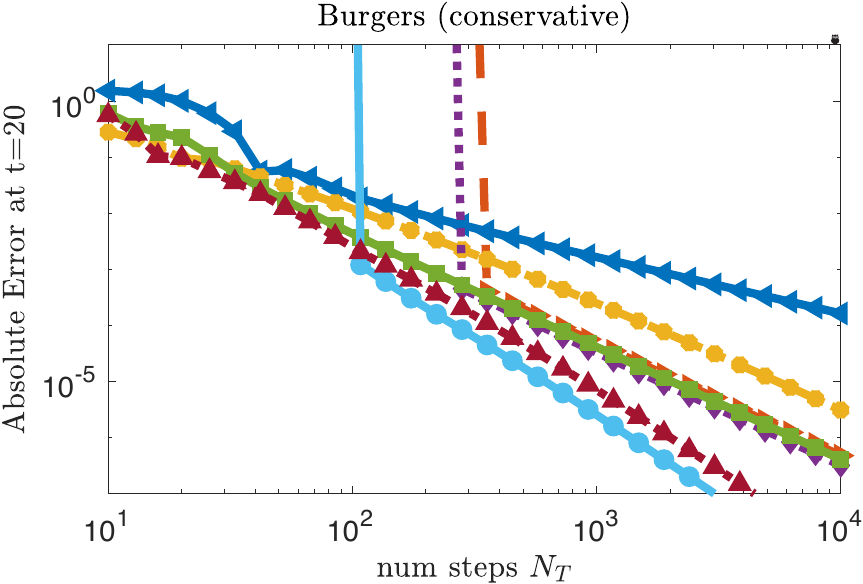}
	
	\caption{Convergence diagrams for the viscous Burgers' equation \cref{eq:burgers} using the nonlinear partitions \cref{eq:burgers-nonlinear-partitioning-non-conservative,eq:burgers-nonlinear-partitioning-conservative} . We discretized the domain $x \in [-8,8]$ using 1000 spatial grid points, and integrate to time $t=20$ using the initial condition $u(x,t=0) = e^{-3x^2}$.}
	\label{fig:burgers-experiment-2}	
\end{figure}

We solve Burgers' equation using each nonlinear partition and compare the convergence of methods that are stable in the coupled stiff $z_2$ limit to those that are not. \Cref{fig:burgers-experiment-2} contains convergence diagrams for four IMEX-NPRK methods that are stable in the coupled stiff $z_2$ limit 
		% first order
		(IMEX-NPRK 1[21] \cref{eq:imex-nprk1-21},
		% second order
		2[42]a \cref{eq:IMEX-NPRK-42},
		2[42]-SiSa \cref{eq:IMEX-NPRK2-43-SiSa},
		% third order
		and 3[54]-Sa \cref{eq:IMEX-NPRK3-54-Sa})
	and three that are not 
		% second order
		(IMEX-NPRK 2[31] \cref{eq:IMEX-NPRK2-31},
		2[42]b \cref{eq:IMEX-NPRK-42},
		% third order
		and 3[54]-Si \cref{eq:IMEX-NPRK354-Si}).
Regardless of the partition, all methods converge at the expected order-of-accuracy; IMEX-NPRK2[31] exhibits excellent error for the non-conservative splitting but remains second-order. For the non-conservative right-hand-side all methods are stable at all timesteps; this is expected since the explicit component of $F(\mathbf{u}, \mathbf{u})$ is truly non-stiff. Conversely, for the conservative right-hand-side we can see that only methods satisfying \cref{eq:stiff-z2-condition} are stable at large timesteps. Moreover, the degree of instability is proportional to the maximum value of the function $\gamma(\theta)$ from \cref{eq:stiff-z2-condition-simplified}. In summary, for problems where the explicit component may contain stiffness, it is important that a method be stable in the coupled, stiff $z_2$ limit.

\subsection{Gray thermal radiation transport}\label{sec:results:trt}

Here we consider an implicit-explicit and semi-implicit formulation of gray thermal radiation transport (TRT). The governing equations for angular intensity $I$ and temperature $T$ are given by:
\begin{subequations}\label{eq:ho}
\begin{align}\label{eq:ho-I}
  \frac{1}{c} \frac{\partial I}{\partial t} & = - \Omega \cdot \nabla I - \sigma_t(T) I + \tfrac{1}{4\pi}ac\sigma_a(T)T^4,\\
  \rho c_v \frac{\partial T}{\partial t} & = - \sigma_p(T) acT^4 + \sigma_E(T)\int I\textnormal{d}\Omega ,\label{eq:ho-T}
\end{align}
\end{subequations}
where $c$ is the speed of light, $c_v$ is the heat capacity, $\rho$ is the material density, and $\sigma(T)$ are opacities. We use a discrete ordinates high-order (HO) low-order (LO) formulation of TRT, with a discontinuous Galerkin discretization of the HO transport equation \eqref{eq:ho} and a finite-volume discretization of the LO moment equations \cite{Park.2020}. Note that the equations are high-dimensional, depending on time, space, and angle due to the collisional integral over direction of transport, $\int I\textnormal{d}\Omega$, and the coupled set of equations is very stiff due to the speed of light scaling. In large-scale simulations such as coupled non-relativistic radiation hydrodynamics, we often want to step over the transport timescale, and as a result \eqref{eq:ho} is treated implicitly in time.

Opacities are nonlinear functions of temperature and density, and can often only be provided via a tabular data lookup. For this reason, it is standard practice in transport simulation to use a first-order semi-implicit integration \cite{Larsen.1988}, where opacities are evaluated on the previous solution value, and all other terms in \eqref{eq:ho} are treated implicitly. Notice that this is exactly a nonlinear partitioning of equations and semi-implicit Euler step as in \eqref{eq:nprk-euler}. By appealing to the NPRK methods developed herein, we can extend this approach to higher order, which we refer to as semi-implicit.
Due to the high dimensionality and stiff behavior, (semi-)implicit integration of \eqref{eq:ho} is very expensive. To alleviate cost, we recently proposed a semi-implicit-explicit approach for integrating \eqref{eq:ho}, based on treating certain asymptotic behavior implicitly and the remaining terms, including opacity evaluations, explicitly \cite{imex-trt}. Each of these approaches for integrating \eqref{eq:ho} are based on a nonlinear partitioning. We now demonstrate the advantages of our new integrators in this context. We compare with the semi-implicit approach of Boscarino \cite{Boscarino.2016,Boscarino.2015} with additive IMEX schemes \emph{LIMEX-Euler} \eqref{eq:nprk-euler}, \emph{H-LDIRK2(2,2,2)} \cite[Table II]{Pareschi.2005}, a two-stage 2nd-order method consisting of 2nd-order L-stable SDIRK implicit and 2nd-order SSP explicit methods, and \emph{SSP-LDIRK2(3,3,2)} \cite[Table IV]{Pareschi.2005}, a three-stage, 2nd-order method consisting of 2nd-order L-stable DIRK implicit and 2nd-order SSP explicit methods. These are the best IMEX-RK schemes that we have found in the literature for these problems \cite{imex-trt,southworth2023implicit}, particularly in terms of stability. We compare with NPRK2[42]b \eqref{eq:IMEX-NPRK-42}, NPRK2[43]-Si \eqref{eq:IMEX-NPRK2-43-Si}, NPRK2[43]-SiSa \eqref{eq:IMEX-NPRK2-43-SiSa}, and NPRK3[54]-Si \eqref{eq:IMEX-NPRK354-Si}. 

\begin{figure}[!htb]
  \centering
    \includegraphics[width=0.6\textwidth]{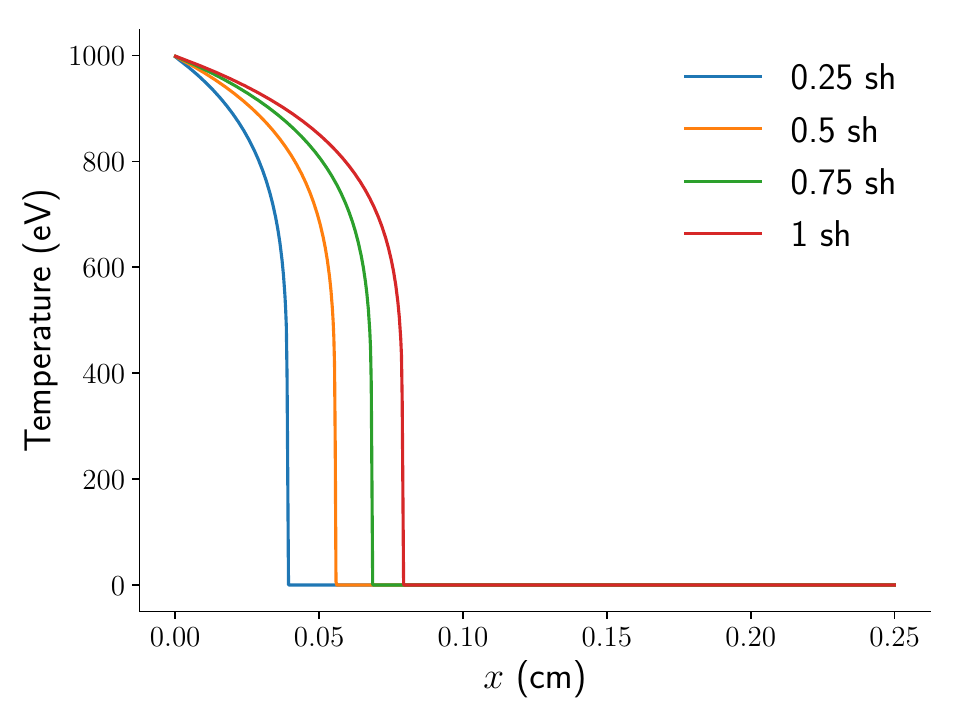}
    \caption{Propagation of Marshak wave from initially uniform background temperature with coarse spatial mesh at four subsequent times in sh (taken from \cite[Fig. 2]{imex-trt}).}
    \label{fig:thick-marshak}
\end{figure}

To test the NPRK methods, we use a classical stiff test problem in TRT, a thick Marshak wave. A Marshak wave is a radiative driven heat wave. Here, the domain is [0,0.25cm], and the radiation and temperature are initially in equilibrium, with a temperature of 0.025 eV. At $x=0$ we have a radiation source corresponding to a Plankian distribution with $T=1000$eV. We use an ideal gas EOS with $\gamma = 5/3$, specific heat $c_V = 3\cdot 10^{12}$, and opacity $\sigma(T) = 10^{12}/T^3$, with initial opacity at $t=0$ $\sigma(0.025)\approx 10^{16}$. As radiation enters the domain, material heats up and the opacity decreases, allowing faster propagation of radiation energy and heat. The evolution of the heat wave at three time snapshots is shown in \Cref{fig:thick-marshak}. We consider two discrete regimes, high spatial resolution ($n=1000$ spatial cells in 1d) and low spatial resolution ($n = 128$ spatial cells in 1d), and both are coupled with an S$_8$ discretization in angle. For brevity, error at the final time 10sh is shown in terms of radiation energy $E = \int I d\Omega$, a quantity of interest due to direct coupling to other physics such as hydrodynamics \cite{southworth2023implicit}; however, analogous convergence results hold for angular intensity $I$ and temperature $T$. For this problem, we simultaneously have to consider order reduction due to stiff nonlinear behavior, order reduction due to partitionings, instability due to partitionings, and instability due to non-physical negativities that arise in temperature. Although negativities can be ``zeroed'' to avoid non-physical or NAN values, a standard approach in practice which we adopt, our experience indicates fixing many negativities can result in loss of convergence or even divergence. 

\begin{figure}[!htb]
  \centering
  \begin{subfigure}[b]{\textwidth}
	  \caption{Fine Spatial Mesh} \rule{\textwidth}{0.4pt} \\[0.5em]
	  \label{subfig:trt-high}
	  \begin{tabular}{cc}
	  		\hspace{2em} {\footnotesize Implicit-explicit} & {\footnotesize Semi-implicit} \hspace{3.5em} \\
	  		\includegraphics[trim={0 0 0 0.82cm},clip,height=1.4in]{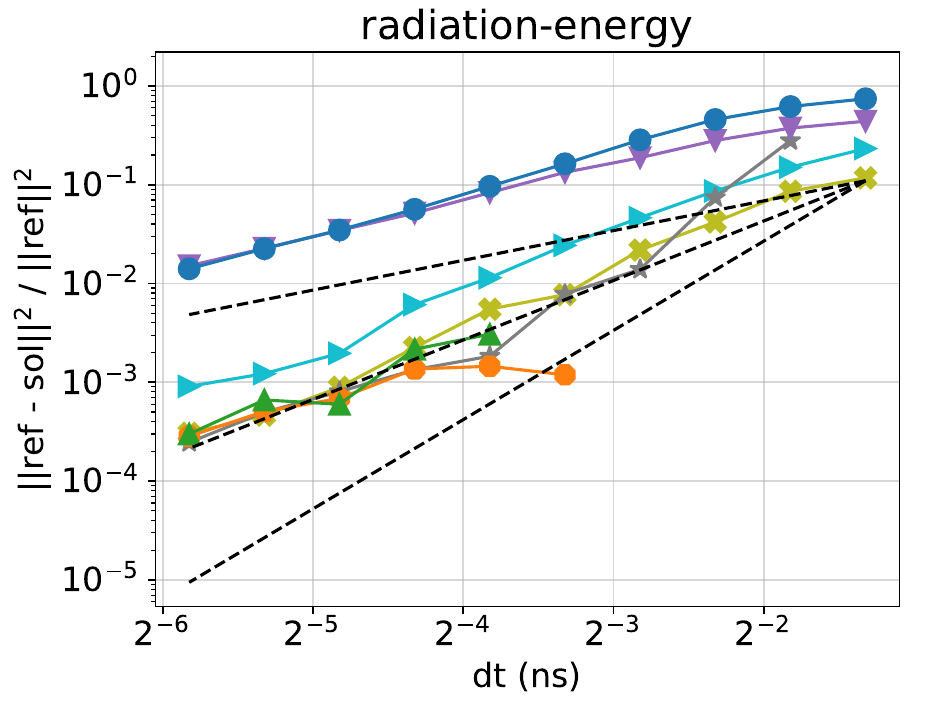} &
	  		\includegraphics[trim={0 0 0 0.90cm},clip,height=1.4in]{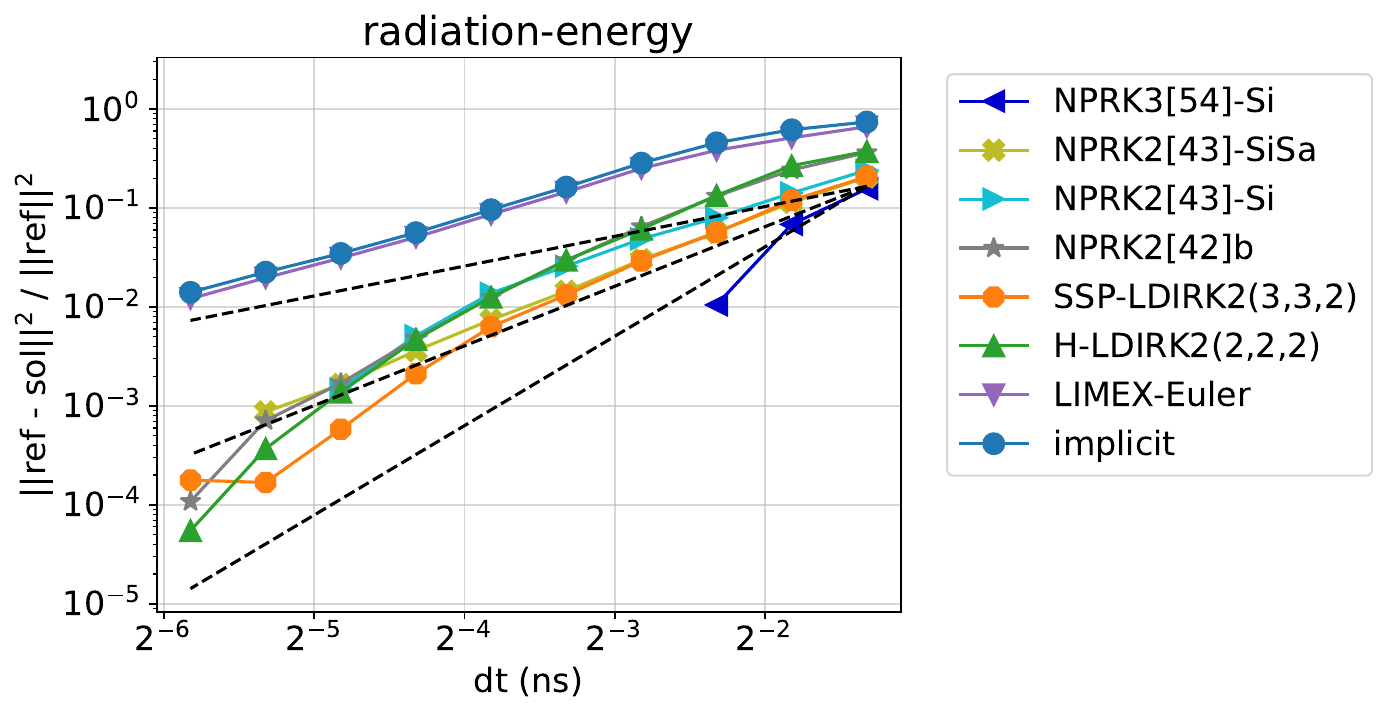}
	  \end{tabular}
 	\end{subfigure}
	\begin{subfigure}[b]{\textwidth}
	  \caption{Coarse Spatial Mesh} \rule{\textwidth}{0.4pt} \\ [0.5em]
	  \label{subfig:trt-low}
	  \begin{tabular}{cc}
	  		\hspace{2em} {\footnotesize Implicit-explicit} & {\footnotesize Semi-implicit} \hspace{3.5em} \\
	  		\includegraphics[trim={0 0 0 0.82cm},clip,height=1.4in]{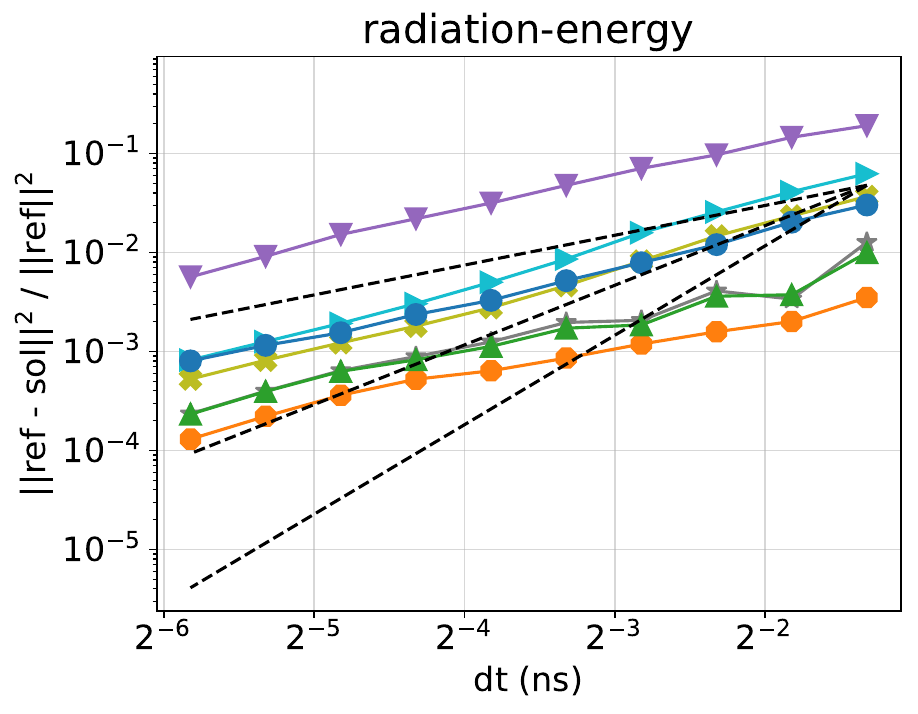} &
	  		\includegraphics[trim={0 0 0 0.90cm},clip,height=1.4in]{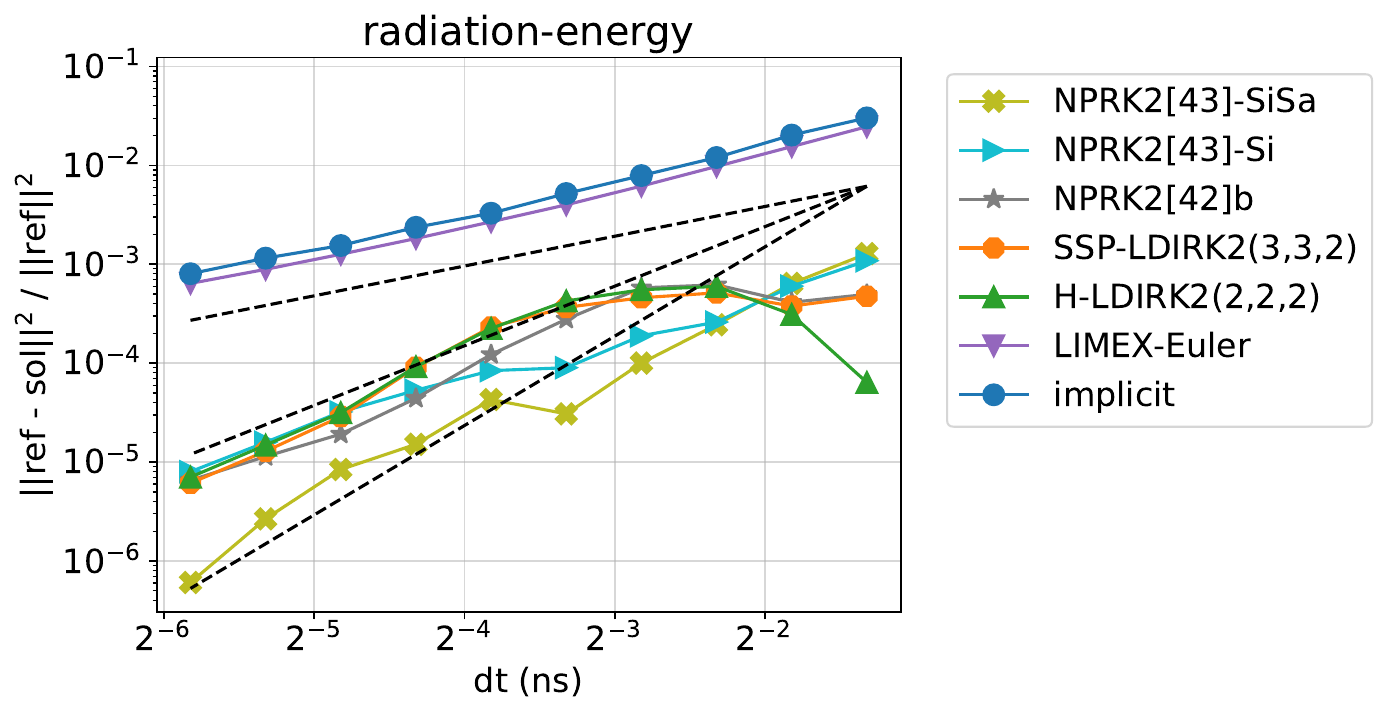}
	  \end{tabular}
 	\end{subfigure}
	
   \caption{Propagation of Marshak wave from initially uniform background temperature solved using a fine spatial grid (top) and coarse spatial grid (bottom). All plots show relative error in radiation energy.}
\end{figure}

High spatial resolution allows us to better resolve the wavefront in space, and thereby its evolution in time; results for high spatial resolution are shown in \Cref{subfig:trt-high}. Overall, we are able to achieve second order accuracy at reasonable physical timesteps for both integration schemes. However, in resolving the wavefront, we introduce an additional challenge of preventing negativities by overstepping the explicit component in temperature, whose discrete stiffness increases with finer spatial resolution (this applies to semi-implicit and IMEX, although is more pronounced in IMEX). As a result, we see that the H-LDIRK and SSP-LDIRK schemes applied to the duplicated set of equations are only stable at smaller timesteps. In contrast, all tested 2nd-order NPRK methods are stable for much larger timesteps, while attaining comparable accuracy at small time steps. Moreover, the three implicit stage NPRK methods require less memory than the three-stage SSP-LDIRK method, due to the design using only implicit stages, and stiffly accurate property for NPRK2[43]-SiSa. The third order scheme NPRK3[54]-Si in most cases develops too many negativities, which results in NaNs and nonphysical solutions. However, for large time steps in the semi-implicit case where it is stable we see nice 3rd order convergence. Note, no other 3rd order schemes we have tested from \cite{Kennedy.2003tv4,Pareschi.2005,Sebastiano.2023} offer stability for these problems either.

% ===============================================
% START - New Combined Figure
% ===============================================

Although sufficient spatial cells to resolve wavefronts is ideal, it is not always feasible in practical TRT simulations, e.g. radiation energy can decay exponentially fast at material interfaces. Results for low spatial resolution are shown in \Cref{subfig:trt-low}. With low spatial resolution, we end up with complete order reduction of the IMEX partitioning to first order, for all tested higher order schemes. This implies to get higher order we must take a semi-implicit approach. In both cases, however, with low spatial resolution, schemes are significantly less sensitive to negativities and instabilities at larger timesteps compared with the high spatial resolution setting. Indeed, all tested schemes converge for all time steps (it should be noted that many schemes of classical IMEX type in Boscarino format are \emph{not} stable). All proposed 2nd order NPRK methods remain robust, with comparable accuracy to the semi-implicit schemes (NPRK offers slightly better accuracy for semi-implicit, but worse for IMEX). NPRK2[43]-SiSa offers the best stability in the sense of basin of convergence, converging at 2nd order even at the largest time steps.

%% file: section-conclusion.tex
\section{Summary, conclusions, and future work}
\label{sec:conclusions}

We introduced a new RK-based framework for solving nonlinearly partitioned equations, analyzed the theoretical properties of the resulting NPRK methods, and constructed IMEX and IMIM NPRK methods with good stability properties when both components possess stiffness. We then demonstrated the efficiency and accuracy of the new methods on the simple one dimensional viscous Burger's equation and on the high dimensional gray thermal radiation transport PIDE. The NPRK framework provides a novel approach for solving ODEs that generalizes classical RK, ARK, and PRK methods. Specifically, the added flexibility for partitioning nonlinear terms allows one to arbitrarily factor nonlinearities and then treat each factor with a different degree of implicitness. We hope that this approach will prove beneficial for a wide range of application problems that are not amenable to additive or component partitioning. 

The class of general NPRK methods \cref{eq:nprk-general-m-partitions} is vast. Future work will investigate the construction of higher-order methods, methods with more varied form of implicitness, and methods with more than two nonlinear components. We also plan to  explore additional application domains that can benefit from nonlinear splittings.

%% file: section-appendix.tex
\begin{appendix}

\section{NPRK sequentially coupled method equivalence}
\label{app:fully-coupled-to-simply-coupled-proof}

\begin{proposition}
	The $s$-stage method \cref{eq:nprk-general-imim-imex} can be expressed as an $\hat{s}$ stage sequentially-coupled method \cref{eq:nprk-simple-coupling-imim-imex} with $\hat{s} \ge s$. 
\end{proposition}

\begin{proof}

	To simplify notation we rename $b_{ij}$ in \cref{eq:nprk-general-imim-imex} to $a_{s+1,i,j}$, so that the output $y_{n+1} = Y_{s+1}$ and all stages can be written as
	\begin{align*}
			Y_i &= y_n + h \left[ a_{i,i,i-1} F(Y_i,Y_{i-1}) + a_{i,i-1,i} F(Y_{i-1},Y_{i}) + \sum_{j=1}^{i-1} \sum_{k=1}^{i-1} a_{ijk} F(Y_j,Y_k) \right].
	\end{align*}
	Similarly, we rename $b_{i,i-1}$ in \cref{eq:nprk-simple-coupling-imim-imex} to $a_{s+1,i,i-1}$.	
	
	We proceed by induction and use tildes to represent variables associated with the sequentially-coupled method. For $s=0$, the statement is trivially satisfied with $\tilde{s}=0$. Next, we assume that the statement holds true for $s=\sigma-2$ and consider the case $s=\sigma-1$. From the assumption, there exist coefficients $\tilde{a}_{i,j,j-1}$ for $i = 2, \ldots, \tilde{s}+1$ so that we can represent stages $Y_1, \ldots, Y_{\sigma-1}$. Therefore, we only need to show that we can find coefficients $\tilde{a}_{i,j,j-1}$, $i > \tilde{s}+1$ to represent the final stage
		\begin{align}
			\begin{split}
			Y_{\sigma}	= y_n + h \bigg[ a_{\sigma,\sigma,\sigma-1} F(Y_{\sigma},Y_{\sigma-1}) + a_{\sigma,\sigma-1,\sigma} F(Y_{\sigma-1},Y_{\sigma}) + \sum_{j=1}^{\sigma-1} \sum_{k=1}^{\sigma-1} a_{\sigma,j,k} F(Y_j,Y_k) \bigg].
			\end{split}
		\end{align}
		The two implicit terms are naturally included in the sequentially coupled ansatz. We therefore focus on the explicit terms and proceed by separating the summed terms into three categories:
		\begin{align}
			k,j < \sigma - 1: & \quad h \sum_{j=1}^{\sigma - 2} \sum_{k=1}^{\sigma - 2} a_{\sigma+1,j,k} F(Y_j,Y_k) && 
				\begin{aligned}
					&\text{\small by induction assumption} \\ 
					& \exists ~ \ell(j,k) \le \tilde{s}+1 \text{ \small so that} \\
					&\text{\small $F(Y_j,Y_k) = F(Y_\ell, Y_{\ell-1})$,}	
				\end{aligned}
			 \\
			j = \sigma - 1: & \quad h \sum_{k=1}^{\sigma - 1} a_{\sigma,\sigma-1,k} F(Y_{\sigma - 1},Y_k) && \text{all terms $F(Y_{\sigma - 1},Y_k)$ are new.} \label{eq:equiv-proof-final-stage-new-terms-I} \\
			k = \sigma - 1: & \quad h \sum_{j=1}^{\sigma-2} a_{\sigma,k,\sigma-1} F(Y_k,Y_{\sigma - 1}) && \text{all terms $F(Y_j,Y_{\sigma - 1})$ are new.} \label{eq:equiv-proof-final-stage-new-terms-II}
		\end{align}
		Note that the sum in \cref{eq:equiv-proof-final-stage-new-terms-II} terminates at $\sigma-2$ instead of $\sigma-1$ to avoid double counting the term $a_{\sigma,\sigma-1,\sigma-1}F(Y_{\sigma-1},Y_{\sigma-1})$ which is already accounted for in \cref{eq:equiv-proof-final-stage-new-terms-I}.	To represent all new terms we append $2\sigma-2$ new stages after $Y_{\tilde{s}+1}$ (i.e. $2\sigma-2$ new rows in $\tilde{a}$)
		\begin{align}
			\begin{aligned}
				&\widetilde{Y}_{\tilde{s} + 2m + 2} = Y_{m+1} \\
				&\widetilde{Y}_{\tilde{s} + 2m + 3} = Y_{\tilde{s}+1} = Y_{\sigma-1}
			\end{aligned} &&
			m = 0, \ldots, \sigma - 2.
		\end{align}
		It follows that
		\begin{align}
			\begin{aligned}
				F(Y_{\kappa},Y_{\sigma-1}) &= F(\widetilde{Y}_\nu,\widetilde{Y}_{\nu-1}) \quad \text{for} \quad \nu =\tilde{s} + 2\kappa, \\
				F(Y_{\sigma-1},Y_{\kappa}) &= F(\widetilde{Y}_\nu,\widetilde{Y}_{\nu-1}) \quad \text{for} \quad \nu =\tilde{s} + 2\kappa+1.
			\end{aligned}
			\label{eq:simply-coupled-additional-stages}
		\end{align}
		These new stages make it possible to express any linear combination $F(Y_{\sigma-1},Y_\kappa)$ or $F(Y_\kappa, Y_{\sigma-1})$ with sequential stages $F(\widetilde{Y}_\nu,\widetilde{Y}_{\nu-1})$. Therefore, we can now add the final stage $\widetilde{Y}_{\tilde{s} + 2\sigma + 1} = Y_\sigma$, which leads a sequentially-coupled method with a total of $(\tilde{s} + 1) + (2\sigma-2) + 1$ stages that is equivalent to the original method.

\end{proof}

\end{appendix}

%% file: section-supplemental.tex
%\appendix

% == TURN APPENDIX INTO SUPPLEMENTAL =============

\def\siamprelabel{SM}
\setcounter{section}{0}
\setcounter{page}{1}

% Does Not Work (without modification \maketitle only prints once)
%\def\siampretitle{Supplementary Materials: }
%\maketitle

% Manual Title
\begin{center}
	{\bfseries{\MakeUppercase{Supplementary Materials: A New Class of Runge-Kutta Methods for Nonlinearly Partitioned Systems}}}
	
	\vspace{1em}
	TOMMASO BUVOLI {\scriptsize AND} BEN S. SOUTHWORTH	
	\vspace{2em}	
\end{center}

% == TURN APPENDIX INTO SUPPLEMENTAL =============

% ====================================================
% APPENDIX: NPRK fully and simply coupled equivalence
% ====================================================

% ====================================================
% APPENDIX: Underlying Methods for Sequentially-Coupled NPRK
% ====================================================

\section{Underlying methods of sequentially-coupled NPRK \cref{eq:nprk-simple-coupling-imim-imex}}
\label{app:sequentially-coupled-underlying-method}

The underlying methods of \cref{eq:nprk-simple-coupling-imim-imex} are
\begin{align}
		a^{\{1\}}_{ij} &= \begin{cases} 0 & j = 1 \\ a_{i,j,j-1} & j > 1 \end{cases} & 
		b^{\{1\}}_i&= \begin{cases} 0 & i = 1 \\ b_{i,i-1} & i > 1 \end{cases} & 
		c^{\{1\}}_i&=\sum_{j=1}^s a^{\{1\}}_{ij} 
		\label{eq:nprk-implicit-method} \\ 
		a^{\{2\}}_{ij} &= \begin{cases} 0 & i = 1 \\ a_{i,j+1,j} & i > 1 \end{cases} &
		b^{\{2\}}_i &= \begin{cases} b_{i+1,i} & i < s \\ 0 & i = s \end{cases}  & 
		c^{\{2\}}_i &=\sum_{j=1}^s a^{\{2\}}_{ij},	
		\label{eq:nprk-explicit-method}
\end{align}
where $a^{\{1\}}$, $a^{\{2\}} \in \mathbb{R}^{s\times s}$ and $b^{\{1\}}$, $b^{\{2\}} \in \mathbb{R}^{s}$. Both methods are reducible since $M^{\{1\}}$ does not use its first stage, and $M^{\{2\}}$ does not use its final stage. The reduced methods are
\begin{align}
		\tilde{a}^{\{1\}}_{ij} &= a_{i+1,j+1,j} & 
		\tilde{b}^{\{1\}}_i&= b_{i+1,i} & 
		\tilde{c}^{\{1\}}_i&=\sum_{j=1}^{s-1} \tilde{a}^{\{1\}}_{ij}  \\ 
		\tilde{a}^{\{2\}}_{ij} &= a_{i,j+1,j} &
		\tilde{b}^{\{2\}}_i &= b_{i+1,i} & 
		\tilde{c}^{\{2\}}_i &=\sum_{j=1}^{s-1} \tilde{a}^{\{2\}}_{ij},
\end{align}
where $\tilde{a}^{\{1\}}, \tilde{a}^{\{2\}} \in \mathbb{R}^{s-1\times s-1}$ and $\tilde{b}^{\{1\}}, \tilde{b}^{\{2\}} \in \mathbb{R}^{s-1}$. As described in \cref{remark:oc-nprk-underlying}, these formula are equivalent to those in \cref{eq:prk-implicit-method,eq:prk-explicit-method}. For clarity we show the tableaux when $s=3$:
\begin{align*}
	&  && \text{\bf Tableau } (a,b,c) && \text{\bf Reduced Tableau } (\tilde{a},\tilde{b},\tilde{c}) \\[1em]	
	& M^{\{1\}}:
	&& \begin{tabular}{r|lll}
		0 						& 0 &  \\
		$(a_{221})$				& 0 & $a_{221}$	\\
		$(a_{221} + a_{332})$ 	& 0 & $a_{321}$ & $a_{332}$	\\ \hline
								& 0 & $b_{21}$  & $b_{32}$
	\end{tabular}
	&& \begin{tabular}{r|ll}
		$(a_{221})$				& $a_{221}$	\\
		$(a_{221} + a_{332})$ 	& $a_{321}$ & $a_{332}$	\\ \hline
								&  $b_{21}$  & $b_{32}$
	\end{tabular}
	\\[1em]
	& M^{\{2\}}:
	&& \begin{tabular}{r|lll}
		& 0 \\
		$(a_{221})$ & $a_{221}$ \\ 
		$(a_{221} + a_{332})$ & $a_{321}$ & $a_{332}$	\\ \hline
		 & $b_{21}$ & $b_{32}$ & $0$
	\end{tabular} 
	&& \begin{tabular}{r|ll}
		& 0 \\
		$(a_{221})$ & $a_{221}$ \\ \hline 
		 & $b_{21}$ & $b_{32}$
	\end{tabular}
\end{align*}

% ====================================================
% APPENDIX: Dahlquist relation to F(Y,Y)
% ====================================================

\section{Reducing nonlinearly partitioned systems to Dahlquist}
\label{app:np-system-dalquist-relation}

Starting from the full nonlinearly partitioned system $y'=F(y,y)$ we locally linearize, neglect nonlinear terms and obtain
\begin{align}
	y' = F(y_n,y_n) + \left[ \frac{\partial F}{\partial y_1}(y_n,y_n) + \frac{\partial F}{\partial y_2}(y_n,y_n) \right] (y - y_n).
\end{align}
Then, letting $y \to y + b$ for $b = y_n - \left[ \frac{\partial F}{\partial y_1}(y_n,y_n) + \frac{\partial F}{\partial y_2}(y_n,y_n) \right]^{-1} F(y_n,y_n)$ yields
\begin{align}
	y' = \left[ \frac{\partial F}{\partial y_1}(y_n,y_n) + \frac{\partial F}{\partial y_2}(y_n,y_n) \right] y
\end{align}
Note that the variable transformation also assumes the matrix $\left(\frac{\partial F}{\partial y_1} + \frac{\partial F}{\partial y_2}\right)$ is invertible. Lastly, if both $\frac{\partial F}{\partial y_1}$ and $\frac{\partial F}{\partial y_2}$ are simultaneously diagonalizable, then we reduce to a decoupled system of partitioned Dahlquist equations \cref{eq:dahlquist-partitioned}.

% ====================================================
% APPENDIX: Proof that B_{\infty}(\epsilon) is polynomial
% ====================================================

\section{Proofs regarding stability functions $B_{\infty}(\epsilon)$ and $\gamma(\theta)$}
\label{app:stability-function-proofs}

	\begin{proposition}
			Consider an $s$-stage diagonally-implicit IMEX ARK method (i.e. $A^{\{1\}}$ lower triangular and $A^{\{2\}}$ strictly lower triangular) where $A^{\{1\}}$ has non-zero diagonals. The function $B_{\infty}(\epsilon) = \lim_{|z_1|\to\infty} R(z_1,\epsilon z_1)$ is a polynomial.	
	\end{proposition}
	\begin{proof}
		The stability function of a diagonally-implicit ARK method is the rational function
		\begin{align}
			R(z_1,z_2) = \frac{p(z_1,z_2)}{q(z_1)} = \frac{\sum_{i=0}^s c_i(z_2) z_1^i}{1 + \sum_{i=1}^{s} d_i z_1^i} \quad \text{for} \quad c_i(z_2) = \left[ \sum_{j=0}^{s-i} c_{ij} z_2^j \right];
			\label{eq:ark-stability-rational-function}
		\end{align}
		see for example \cite{Kennedy.2003tv4}. The coefficients $d_i$, $c_{ij}$ are real numbers that depend on the method coefficients, with the denominator	having the simple form
		\begin{align}
			q(z) = \prod_{i=1}^s (1-A^I_{ii} z).	
		\end{align}
		To derive $B_\infty(\epsilon)$ we first substitute $z_2 = \epsilon z_1$ in the numerator of $R(z_1,z_2)$
		\begin{align}
			\text{Numerator}[R(z_1,\epsilon z_1)] =  p(z_1,\epsilon z_1) = \sum_{i=0}^s \left[ \sum_{j=0}^{s-i} c_{ij} \epsilon^j z_1^j \right] z_1^i.
		\end{align}
		Now rewrite the sum in standard powers of $z_1$,
		\begin{align}
			\sum_{i=0}^s \left[ \sum_{j=0}^{s-i} c_{ij} \epsilon^j z_1^j \right] z_1^i = \sum_{i=0}^{s} \sum_{p=i}^s c_{i,p-i} \epsilon^{p-i} z_1^p, \quad p = i+j.
		\end{align}
		Next, using $\sum_{i=0}^s \sum_{p=i}^s a_{i,p} = \sum_{p=0}^s \sum_{i=0}^p a_{i,p}$ we have,
		\begin{align}
			\sum_{i=0}^s \sum_{p=i}^{s} c_{i,p-i} \epsilon^{p-i} z_1^p = \sum_{p=0}^s z_1^p \left[ \sum_{i=0}^p c_{i,p-i} \epsilon^{p-i} \right] = \sum_{p=0}^s \tilde{c}_p(\epsilon) z_1^p.
		\end{align}
		Lastly, we consider the denominator
		\begin{align}
			\text{Denomenator}[R(z_1,\epsilon z_1)]	= \prod_{i=1}^s (1-A^I_{ii} z_1) = 1 + \sum_{i=1}^{s} d_i z_1^i.
		\end{align}
		Using $d_s = (-1)^s \prod_{i=1}^s A^{\{1\}}_{ii}\ne 0$ (since $A^{\{1\}}_{ii} \ne 0$ by assumption), we have 
		\begin{align}
			B_{\infty}(\epsilon) = \lim_{|z_1|\to \infty} R(z_1,\epsilon z_1) = \tilde{c}_s(\epsilon) / d_s,
		\end{align}
		which is a polynomial in $\epsilon$.
	\end{proof}
	
	\begin{proposition}
			The function $\gamma(\theta) = |\beta_{\infty}(e^{i\theta})|^2$ is a cosine series.
	\end{proposition}
	\begin{proof}
	This follows from the fact that 	the squared magnitude of any polynomial with real-valued coefficients evaluated on the unit circle is a cosine series. The following algebraic proof also reveals the coefficients. Let $p(z) = \sum_{j=0}^s c_j z^j$ for $c_j \in \mathbb{R}$, then $|p(z)|^2 = p(z)(p(z))^* = p(z)p(z^*)$. Next,
	\begin{align}
		 \left( \sum_{j=0}^s c_j e^{i j \theta} \right)
		 \left( \sum_{k=0}^s c_k e^{-i k \theta} \right) 
		 &= \sum_{j=0}^s \sum_{k=0}^s c_j c_k e^{i(j-k) \theta} \\
		 &= \sum_{j=0}^s \sum_{n=j-s}^j c_j c_{j-n} e^{i n\theta} \quad n = j - k \\
		 &= \sum_{j=0}^{s} c_j c_{j} + \sum_{n=1}^s \left[\sum_{j=n}^{s}  c_j c_{j-n} \left( e^{in\theta} + e^{-in\theta}\right) \right].
	\end{align}
	Lastly, using $\cos(x) = (e^{ix} + e^{-ix})/2$, we have
	\begin{align}
		 \sum_{n=0}^s d_n \cos(n \theta) \quad \text{for} \quad
		 d_n  = 
		 \begin{cases}
		 	\sum_{j=0}^{s}c_{j}^2 & n = 0, \\
		 	\sum_{j=n}^{s} 2 c_j c_{j-n}	 & n > 0	.
		 \end{cases}		
	\end{align}
	\end{proof}

\section{Sequentially-coupled IMEX-NPRK method derivations}
\label{app:method-derivations-imex}

This section contains detailed method derivations for a IMEX sequential coupled methods \cref{eq:nprk-imex-simply-coupled-with-tableau} whose ansatz and tableau are copied below for convenience.
\begin{align*}
	\begin{small}
		\begin{aligned}
			Y_i &= y_n + h\sum_{j=2}^i a_{i,j,j-1} F(Y_j,Y_{j-1}), \quad i=1,\ldots, s, \\
			y_{n+1} &= y_n + h \sum_{j=2}^s b_{j,j,j-1} F(Y_j,Y_k).
		\end{aligned}
	\end{small} && 
	\begin{small}
		\begin{tabular}{l|llll}
			& 0 \\
			& $a_{221}$ & \\
			& $\vdots$ & $\ddots$ & \\ 
			& $a_{s21}$ & $\hdots$ & $a_{s,s,s-1}$ \\ \hline
			& $b_{21}$  & $\hdots$ & $b_{s,s-1}$
		\end{tabular}
	\end{small}
\end{align*}

\subsection{Second-order methods}

Sequentially-coupled IMEX NPRK methods with three and four stages can achieve at most second-order accuracy and require either one, two, or three implicit solves. In our derivations, we will often consider the third-order residual
\begin{align}
	\begin{scriptsize}
	\mathbf{r}^{(3)} = \left[ 
	\renewcommand*{\arraystretch}{1.5}
	\begin{aligned}
		&\textstyle \sum_{i=1}^{s} b_i \left( c_{i} \right)^2 - \tfrac{1}{3} \\
		&\textstyle \sum_{i=1}^{s} b_i c_{i} \widehat{c}_{i} -\tfrac{1}{3} \\
		&\textstyle \sum_{i=1}^{s} b_i \left( \widehat{c}_{i} \right)^2 -\tfrac{1}{3} \\
		&\textstyle \sum_{i=1}^{s} \sum_{j=1}^{s} b_i a_{ij} c_{j} - \tfrac{1}{6} \\
		&\textstyle \sum_{i=1}^{s} \sum_{j=1}^{s} b_i a_{ij} \widehat{c}_{j} - \tfrac{1}{6} \\
		&\textstyle \sum_{j=1}^{s} b_i \widehat{a}_{ij} c_{j} - \tfrac{1}{6} \\
		&\textstyle \sum_{j=1}^{s} b_i \widehat{a}_{ij} \widehat{c}_{j} - \tfrac{1}{6}
	\end{aligned}
	\right],
	\end{scriptsize}
	\label{eq:third-order-residual}
\end{align}
that measures how well a method satisfies the third-order conditions from \cref{tab:PRK-order-conditions}.

\subsubsection{Three-stages, one implicit solve}
\label{app:method-imp-midpoint-exp-midpoint}

The starting ansatz is a three stage, sequentially coupled IMEX-NPRK method that only requires a single implicit solve. The corresponding tableau is
\begin{align}
	\begin{tabular}{l|lll}
		& 0 \\
		& $a_{221}$ & \\
		& $a_{321}$ & 0 \\ \hline
		& $b_{21}$  & $b_{32}$
	\end{tabular}
\end{align}
Order conditions for a second-order method require that
\begin{align}
	b_{21} + b_{32} = 1, && 
	b_{32} a_{221} = \tfrac{1}{2}, && 
	b_{21} a_{221} + b_{32} a_{321}= \tfrac{1}{2}.
\end{align}
This leads to a one-parameter family of methods
\begin{align}
	a_{221} &= \frac{1}{2 a_{b_{32}}}, &
	a_{321} &= \frac{2 b_{32}-1}{2 b_{32}^2}, &
	b_{21} &= 1-b_{32},
\end{align}
with stability function
\begin{align}
	R(z_1,z_2) = \frac{ \left(b_{23} z_2^2+2 b_{23} z_2+2 b_{23}\right) +  \left( 2 b_{23} (z_2+1)-z_2-1 \right)z_1 + \left(b_{23}-1\right)z_1^2}{2 b_{23}-z_1}.
\end{align}
To prevent the stability function from growing unbounded as $|z_1| \to \infty$ we require that the $z_1^2$ coefficient in the numerator be zero
\begin{align}
	(b_{23} - 1) = 0 
	\quad \iff \quad
	b_{23} = 1.
\end{align}
This leads to a method 
\begin{align}
	\begin{tabular}{l|lll}
		& 0 \\
		& $1/2$ & \\
		& $1/2$ & 0 \\ \hline
		& $0$  & $1$
	\end{tabular}
	\label{eq:app-method-imp-midpoint-exp-midpoint}
\end{align}
that is \underline{conditionally} $A$-stable in the stiff $z_1$ limit, specifically
\begin{align}
\lim_{|z_1| \to \infty} R(z_1,z_2) = -(1 + z_2).
\end{align}
Moreover, the underlying methods of \cref{eq:app-method-imp-midpoint-exp-midpoint} are the implicit and explicit midpoint methods. The third-order residual \cref{eq:third-order-residual} for this method and its two norm are
\begin{align}
	\mathbf{r}^{(3)} &= \left\{-\frac{1}{12},-\frac{1}{12},-\frac{1}{12},\frac{1}{12},-\frac{1}{6},\frac{1}{12},-\frac{1}{6}\right\},	\\
	\| \mathbf{r}^{(3)} \|_2 &= \frac{\sqrt{13}}{12} \approx 0.300463.
\end{align}

\subsubsection{Three-stages, two implicit solves}
\begin{figure}[h]
	\begin{subfigure}[b]{0.45\textwidth}
		\includegraphics[width=\textwidth]{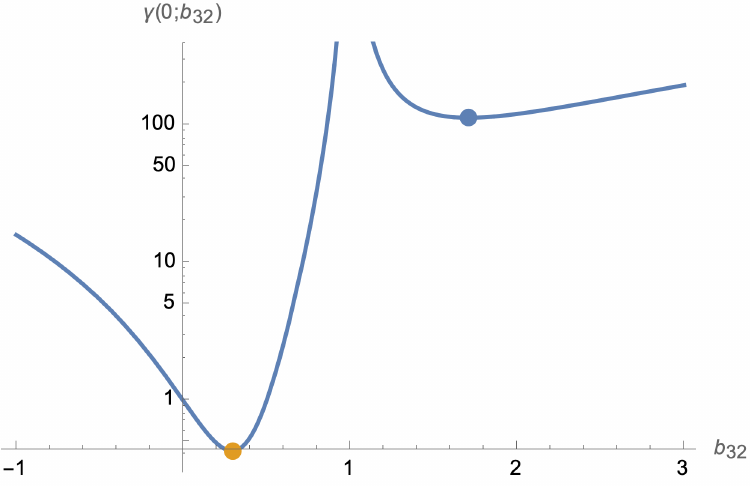}
		\caption{ Function $\gamma(0;b_{32})$ -- log $y$ axis. }
	\end{subfigure}
	\hfill
	\begin{subfigure}[b]{0.45\textwidth}
		\includegraphics[width=\textwidth]{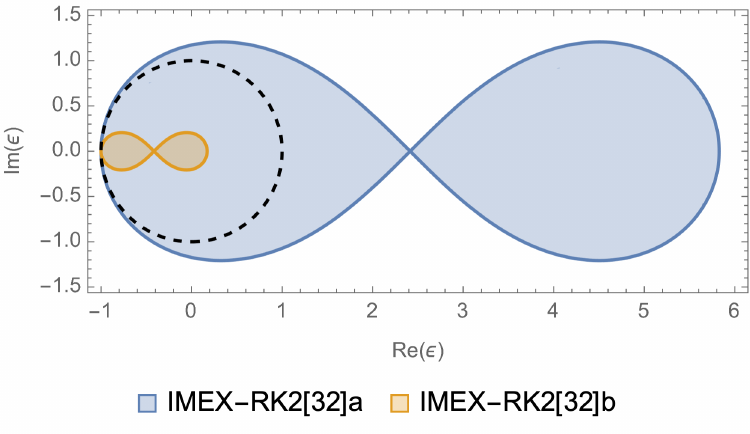}	
		\caption{ Regions $|\beta_\infty(\epsilon) | \le 1$}
	\end{subfigure}

	\caption{Stability figures for the IMEX-NPRK2[32]a and IMEX-NPRK2[32]b methods. (a) shows the function $\gamma(\theta;b_{32})$, defined in \cref{eq:imex-232-gamma}; the local min (blue circle) is $\gamma\left(0; \tfrac{2-\sqrt{2}}{2}\right) = 57+40 \sqrt{2}$, while the global min (orange circle) is  $\gamma\left(0; \tfrac{2+\sqrt{2}}{2}\right) = 57-40 \sqrt{2}$ (b) Stability regions in the stiff $z_1$ limit for IMEX-NPRK2[32] methods with $b_{32} = 1\pm\frac{1}{\sqrt{2}}$.}
	\label{fig:imex-232-stability-extra}
	
\end{figure}

The starting ansatz is a three stage, sequentially coupled IMEX-NPRK method requiring two implicit solves. The corresponding tableau is
\begin{align}
	\begin{tabular}{l|lll}
		& 0 \\
		& $a_{221}$ & \\
		& $a_{321}$ & $a_{332}$ & \\ \hline
		& $b_{21}$  & $b_{32}$
	\end{tabular}
\end{align}
Order conditions for a second-order method require that
\begin{align}
	b_{21} + b_{32} = 1, && 
	b_{21} a_{221} = \tfrac{1}{2}, && 
	b_{21} a_{221} + b_{32} (a_{332}+a_{321}) = \tfrac{1}{2}.
\end{align}
This leads to a two-parameter family of methods
\begin{align}
	a_{221} &= \frac{1}{2 b_{32}}, &
	a_{332} &= \frac{-2 a_{321} b_{32}^2+2 b_{32}-1}{2 b_{32}^2}, &
	b_{21} &= 1-b_{32}.	
\end{align}
The resulting method is not L-stable in the stiff $z_1$ limit, instead satisfying
\begin{align}
	\lim_{|z_1|\to \infty} R(z_1,z_2) = \frac{4 a_{321} b_{32}^3-2 a_{321} b_{32}^2+2 b_{32}^3-6 b_{32}^2+4 b_{32}-1}{-2 a_{321} b_{32}^2+2 b_{32}-1}	.
\end{align}
Enforcing L-stability leads to the one-parameter family of methods
\begin{align}
	a_{221} &= \frac{1}{2 b_{32}}, &
	a_{321} &= \frac{-2 b_{32}^3+6 b_{32}^2-4 b_{32}+1}{2 b_{32}^2 (2 b_{32}-1)}, \\
	a_{332} &= \frac{b_{32}-1}{2 b_{32}-1}, &
	b_{21}  &= 1-b_{32}.
\end{align}
The free parameter $b_{32}$ is selected to minimize $\gamma(\theta; b_{32})$, defined in \cref{eq:stiff-z2-condition-simplified}, on the interval $\theta\in[0,2\pi)$. For the L-stable family of methods we have
\begin{align}
	\gamma(\theta; b_{32}) &= c_0(b_{32}) + c_1(b_{32}) \cos(\theta) \label{eq:imex-232-gamma}, \\ 
	c_0(b_{32}) &= \frac{b_{32} (b_{32} (4 b_{32} (2 b_{32}-3)+9)-4)+1}{(b_{32}-1)^2} , \\
	c_1(b_{32}) &= \frac{2 b_{32} (2 b_{32} (b_{32} (2 b_{32}-3)+2)-1)}{(b_{32}-1)^2}.	
\end{align}
Differentiating $\gamma(\theta; b_{32})$ with respect to $\theta$ reveals that the maximum value occurs at either $\theta = 0$ or $\theta = \pi$. It follows that $\gamma(\pi;b_{32}) = 1$ for all $b_{32}$, so stability at this point cannot be improved. We can however optimize stability at $\theta = 0$ by selecting $b_{32}$ so that
\begin{align}
	\frac{\partial \gamma(0; b_{32})}{\partial b_{32}} = 0,
\end{align}
which leads to the solutions
\begin{align}
	\text{global min: } & b_{32} = 1-\frac{1}{\sqrt{2}},  &
	\text{local min: } & b_{32} = 1+\frac{1}{\sqrt{2}}.
\end{align}
A graph showing $\gamma(\pi;b_{32})$ and the two minimums is contained in \cref{fig:imex-232-stability-extra}. The two-norms of the third-order residuals for the optimal methods are:
\begin{align}
	\label{eq:nprk232s-a-error}
	\|\mathbf{r}^{(3)}(b_{32} = \tfrac{2 - \sqrt{2}}{2})\|_2 &= \frac{1}{6} \left(70 \sqrt{2}+99\right) \sqrt{64433-45561 \sqrt{2}} \approx 4.15904, \\
	\label{eq:nprk232s-b-error}
	\|\mathbf{r}^{(3)}(b_{32} = \tfrac{2 + \sqrt{2}}{2})\|_2 &= \frac{1}{6} \left(99-70 \sqrt{2}\right) \sqrt{64433+45561 \sqrt{2}} \approx 0.302179.
\end{align}

\subsubsection{Four-stages, two implicit solves}

The starting ansatz is a four-stage, sequentially coupled, IMEX-NPRK method. To reduce free parameters we consider a simplified anstaz with $a_{i,3,2}=0$ and tableau
\begin{align}
	\begin{tabular}{l|lll}
		& 0 \\
		& $a_{221}$ & \\
		& $a_{321}$ & $0$ & \\ 
		& $a_{421}$ & $0$ & $a_{443}$ \\ \hline
		& $b_{21}$  & $0$ & $b_{43}$
	\end{tabular}.
\end{align}
Order conditions for a second-order method require that
\begin{align}
	b_{21} + b_{43} &= 1 ,&
	a_{321} b_{43} &= \tfrac{1}{2} ,&	
	b_{21} a_{221} + b_{43} ( a_{421} + a_{443}) &= \frac{1}{2}.
\end{align}
This leads to a three-parameter family of methods
\begin{align}
	a_{321} &= \frac{1}{2 b_{43}}, &
	b_{21}  &= 1-b_{43}, &
	a_{443} &= \frac{2 a_{221} (b_{43}-1)-2 a_{421} b_{43}+1}{2 b_{43}}.	
\end{align}
The resulting method is not L-stable in the stiff $z_1$ limit, instead satisfying
\begin{align}
\lim_{|z_1|\to \infty} R(z_1,z_2) = \frac{a_{221} a_{443}-a_{221} b_{43}+a_{421} b_{43}-a_{443} b_{21}}{a_{221} a_{443}}.
\end{align}
Enforcing L-stability leads to a two-parameter family of methods
\begin{align}
	\begin{aligned}
		a_{421} &= \frac{a_{221} (2 a_{221} (b_{43}-1)-4 b_{43}+3)+b_{43}-1}{2 (a_{221}-1) b_{43}}, \quad
		a_{432} = 0, \\
		a_{443} &= \frac{1}{2 (a_{221}-1)}+1, \quad
		b_{21} = 1-b_{43}.
	\end{aligned}
\end{align}
There are now two free parameters $b_{43}$ and $a_{221}$. We proceed by attempting to minimize $\gamma(\theta; a_{221})$ on $\theta \in [0, 2\pi]$ which only depends on $a_{221}$. Specifically, we have
\begin{align}
	\gamma(\theta; a_{221}) &= c_0(a_{221}) + c_1(a_{221}) \cos(\theta), \label{eq:imex-242-gamma} \\ 
	c_0(a_{221}) &= \frac{a_{221} (a_{221} (4 (a_{221}-2) a_{221}+9)-6)+2}{(1-2 a_{221})^2 a_{221}^2}, \\
	c_1(a_{221}) &= \frac{2-2 a_{221} (2 (a_{221}-2) a_{221}+3)}{(1-2 a_{221})^2 a_{221}^2}	.
\end{align} 
Differentiating $\gamma(\theta; a_{221})$ with respect to $\theta$ reveals that the maximum value occurs at either $\theta = 0$ or $\theta = \pi$. It follows that $\gamma(\pi;a_{221}) = 1$ for all $a_{221}$, so stability at this point cannot be improved. We can however optimize stability at $\theta = 0$ by selecting $a_{221}$ so that
\begin{align}
	\frac{\partial \gamma(0; a_{221})}{\partial a_{221}} = 0.
\end{align}
which leads to the solutions
\begin{align}
	\text{global min: } & a_{221} = 1+\frac{1}{\sqrt{2}}  &
	\text{local min: } & a_{221} = 1-\frac{1}{\sqrt{2}} 
\end{align}
The remaining variable $b_{43}$ is selected to minimize the two norm of the third-order residuals \cref{eq:third-order-residual}. Specifically we select $b_{43}$ which satisfies
\begin{align}
	\frac{\partial \mathbf{r}^{(3)}(b_{43};a_{221})}{\partial b_{43}} = 0.
\end{align}
This leads to the methods
\begin{align}
	\begin{tabular}{ccc}
		\text{\underline{IMEX-NPRK2[42]a}} & \hspace{3em} & \text{\underline{IMEX-NPRK2[42]b}}	\\[0.5em]
		\text{(Global Optimum)} & & \text{(Local Optimum)} \\[0.5em]	
		\renewcommand\arraystretch{1.75}
		\begin{tabular}{l|llll}
			& $0$ \\
			& $1+\frac{1}{\sqrt{2}}$		\\
			& $\frac{26 - 3 \sqrt{2}}{42}$	& $0$ \\
			& $\frac{-20 -23 \sqrt{2}}{42}$ & $0$ & $1+\frac{1}{\sqrt{2}}$ \\ \hline 
			& $\frac{16-9 \sqrt{2}}{94}$ 		& $0$ & $\frac{78+9 \sqrt{2}}{94}$ \\ 
		\end{tabular}
		& & \renewcommand\arraystretch{1.75}
		\begin{tabular}{l|llll}
			& $0$ \\
			& $1-\frac{1}{\sqrt{2}}$		\\
			& $\frac{26 + 3 \sqrt{2}}{42}$	& $0$ \\
			& $\frac{-20 + 23 \sqrt{2}}{42}$ & $0$ & $1-\frac{1}{\sqrt{2}}$ \\ \hline 
			& $\frac{16 + 9 \sqrt{2}}{94}$ 		& $0$ & $\frac{78-9 \sqrt{2}}{94}$ \\ 
		\end{tabular}
	\end{tabular}
\end{align}
The two-norm of the third-order residuals for these optimal methods are:
\begin{align}
	\label{eq:nprk242s-a-error}
	\|\mathbf{r}^{(3)}(a_{221} = \tfrac{2 + \sqrt{2}}{2}, b_{43}=\tfrac{78+9\sqrt{2}}{94})\|_2 &= \frac{1}{6} \sqrt{\frac{1}{7} \left(367 + 253 \sqrt{2}\right)} \approx 1.69593, \\
	\label{eq:nprk242s-b-error}
	\|\mathbf{r}^{(3)}(a_{221} = \tfrac{2 - \sqrt{2}}{2}, b_{43}=\tfrac{78-9\sqrt{2}}{94})\|_2 &= \frac{1}{6} \sqrt{\frac{1}{7} \left(367-253 \sqrt{2}\right)} \approx 0.191112.
\end{align}
which is 40\% and 63\% of the residuals for the respective three-stage two-solve methods with errors \cref{eq:nprk232s-a-error,eq:nprk232s-b-error}.

\subsubsection{Four-stages, three implicit solves}

The starting ansatz is a four-stage, sequentially coupled, IMEX-NPRK method with Tableau
\begin{align}
	\begin{tabular}{l|lll}
		& 0 \\
		& $a_{221}$ & \\
		& $a_{321}$ & $a_{332}$ & \\ 
		& $a_{421}$ & $a_{432}$ & $a_{443}$ \\ \hline
		& $b_{21}$  & $b_{32}$ & $b_{43}$
	\end{tabular}.
\end{align}
Imposing conditions for a second-order method leads to a six-parameter family
\begin{align}
	\begin{scriptsize}
		\begin{aligned}
				a_{332} &= \frac{-2 a_{221} b_{32}-2 a_{321} b_{43}+1}{2 b_{43}}, \quad
			b_{21} = -b_{32}-b_{43}+1 \\
			a_{443} &= \frac{2 a_{221} b_{32}^2+2 a_{221} b_{32} b_{43}+2 a_{221} b_{43}^2-2 a_{221} b_{43}-2 a_{421} b_{43}^2-2 a_{432} b_{43}^2-b_{32}+b_{43}}{2 b_{43}^2}.
		\end{aligned}
	\end{scriptsize}
	\label{eq:four-stage-three-solve-second-order-family}
\end{align}
The resulting method will be L-stable in $z_1$ if we impose the additional condition
\begin{align}
	\begin{split}
		& \frac{1}{a_{221} a_{332} a_{443}} \bigg[ a_{332} (a_{443}-b_{43})+a_{432} b_{43}-a_{443} b_{32}) -a_{321} a_{432} b_{43} \\
		& \hspace{6em} +a_{321} a_{443} b_{32}+a_{332} a_{421} b_{43} -a_{332} a_{443} b_{21} \bigg] = 0.
	\end{split}	
	\label{eq:l-stable-three-solve-second-order-condition}
\end{align}
We are unable to get close form coefficients using Mathematica if we attempt to impose \cref{eq:l-stable-three-solve-second-order-condition}; therefore we seek further simplifications. If we force a SDIRK method ($a_{221}=a_{332}=a_{443}=\gamma$), and then impose L-stability in $z_1$, we obtain the three parameter class of methods
	\begin{align}
		\begin{small}
		\begin{aligned}
				a_{321} &= \frac{1-2 \gamma  (b_{32}+b_{43})}{2 b_{43}}, \quad
			a_{432} = \frac{\gamma  (-2 (\gamma -2) \gamma -1)}{2 \gamma  (b_{32}+b_{43})-1}, \quad
			b_{21} = -b_{32}-b_{43}+1, \\
			a_{421} &= \frac{1}{2} \left(\frac{b_{32} (2 b_{32} \gamma -1)}{b_{43}^2}+\frac{2 (b_{32}-1) \gamma +1}{b_{43}}+\frac{2 \gamma  (2 (\gamma -2) \gamma +1)}{2 \gamma  (b_{32}+b_{43})-1}\right).
		\end{aligned}
		\end{small}
	\end{align}
To derive the optimized method \cref{eq:IMEX-NPRK2-43-Si} we numerically solve for  the remaining three parameters that minimize the objective function
\begin{align}
	c(\gamma,b_{32},b_{43}) = \|\mathbf{r}^{3}(\gamma,b_{32},b_{43})\|^2_2 + \int_{0}^{2\pi} \gamma(\theta; \gamma,b_{32},b_{43}) d\theta.
\end{align}
The first term measures the third-order residual from \cref{eq:third-order-residual}, while the second term measures the area under the always positive stability function $\gamma(\theta)$ from  \cref{eq:stiff-z2-condition-simplified} (ideally we want $\gamma(\theta) \le 1$ for all $\theta \in [0, 2\pi)$ -- minimizing area is a heuristic for this condition that is also simple to optimize over).

An alternative approach for deriving four-stage, three solve methods is to restrict to a stiffly accurate methods (i.e. $b_{i,i-1} = a_{4,i,i-1}$, $i = 2, \ldots, 4$). This leads to the three parameter family
\begin{align}
	\begin{scriptsize}
		\begin{aligned}
				a_{221} &= \frac{a_{432}+2 a_{443}^2-a_{443}}{2 \left(a_{432}^2+a_{432} a_{443}+a_{443}^2-a_{443}\right)}, \quad
				a_{421} = -a_{432}-a_{443}+1 \\
				a_{332} &= \frac{-2 a_{321} a_{432}^2-2 a_{321} a_{432} a_{443}-2 a_{321} a_{443}^2+2 a_{321} a_{443}-2 a_{432} a_{443}+2 a_{432}+a_{443}-1}{2 \left(a_{432}^2+a_{432} a_{443}+a_{443}^2-a_{443}\right)}.
		\end{aligned}
	\end{scriptsize}
	\label{eq:four-stage-three-solve-second-order-stiff-family}
\end{align}
This family of methods is naturally L-stable in $z_1$ and all have stability function $\beta_\infty(\epsilon) = -\epsilon^3$ implying that they are also stable in the coupled stiff $z_2$ limit. If we further impose SDIRK method ($a_{221}=a_{332}=a_{443}=\gamma$) then we obtain two one-parameter method families. These two method families are
\begin{align}
	\begin{small}
		\begin{aligned}
			a_{321} &= \frac{1 -2 \gamma^2 +  f(\gamma)}{4\gamma}, &	
			a_{421} &= \frac{-1 + 4\gamma -2 \gamma^2 + f(\gamma)}{4\gamma},	&
			a_{432} &= \frac{1 - 2 \gamma^2 -  f(\gamma)}{4\gamma}
		\end{aligned}
	\end{small}
	\label{eq:imex-nprk-43-sa-family1} \\
	\begin{small}
		\begin{aligned}
			a_{321} &= \frac{1 -2 \gamma^2 -  f(\gamma)}{4\gamma} &	
			a_{421} &= \frac{-1 + 4\gamma -2 \gamma^2 - f(\gamma)}{4\gamma}	&
			a_{432} &= \frac{1 - 2 \gamma^2 +  f(\gamma)}{4\gamma},
		\end{aligned}
	\end{small}
	\label{eq:imex-nprk-43-sa-family2}
\end{align}
with $f(\gamma) = \sqrt{1-4 \gamma^2 (\gamma (3 \gamma-8)+3)}$. To derive the method \cref{eq:IMEX-NPRK2-43-SiSa} we considered \cref{eq:imex-nprk-43-sa-family1} and numerically solved for the $\gamma$ value that minimizes the objective function
\begin{align}
	c(\gamma) = \|\mathbf{r}^{3}(\gamma)\|_2.
	\label{eq:objective-nprk-43-sa}
\end{align}
This leads to $\gamma = 0.386585$ and the method
\begin{align}
	\begin{tabular}{c}
		\text{\underline{IMEX-NPRK2[43]a-SiSa}}	\\[0.5em]	
		\begin{small}
		\renewcommand\arraystretch{1.25}
		\begin{tabular}{l|llll}
			&			   	   0   &                   &                    \\
			&   0.386585000000000  &                   &                    \\
			&   1.027233588987035  & 0.386585000000000 &                    \\
			&   0.733856970649542  &-0.120441970649542 &  0.386585000000000 \\ \hline 
			&   0.733856970649542  &-0.120441970649542 &  0.386585000000000
		\end{tabular}
		\end{small}
	\end{tabular}
\end{align}
which has $c(386585)=0.500262$. Minimizing \cref{eq:objective-nprk-43-sa} for \cref{eq:imex-nprk-43-sa-family2}, we obtain $\gamma = 0.325754$, and the method
\begin{align}
	\begin{tabular}{c}
		\text{\underline{IMEX-NPRK2[43]b-SiSa}}	\\[0.5em]	
		\begin{small}
		\renewcommand\arraystretch{1.25}
		\begin{tabular}{l|llll}
			&			   	   0   &                   &                    \\
			&   0.325754000000000  &                   &                    \\
			&   -0.036442462428244  & 0.325754000000000 &                    \\
			&   -0.571343031569375  &1.245589031569375 &  0.325754000000000 \\ \hline 
			&   -0.571343031569375  &1.245589031569375 &  0.325754000000000
		\end{tabular}
		\end{small}
	\end{tabular}
	\label{eq:IMEX-NPRK2[43]b-SDIRK-SA}	
\end{align}
which has $c(0.325754)=0.286004$.

\subsection{Third-order methods}

We require at least five stages to derive a third-order method using the sequentially coupled IMEX NPRK ansatz \cref{eq:nprk-imex-simply-coupled-with-tableau}.

\subsubsection{Five-stages, four implicit solves}

	We consider a method with tabaleau
	\begin{align}
		\begin{tabular}{l|llll}
		& 0 \\
			& $a_{221}$ & \\
			& $a_{321}$ & $a_{332}$ & \\ 
			& $a_{421}$ & $a_{432}$ & $a_{443}$ \\ 
			& $a_{521}$ & $a_{532}$ & $a_{543}$ & $a_{554}$ \\ \hline
			& $b_{21}$  & $b_{32}$ & $b_{43}$ & $b_{54}$
		\end{tabular}
	\end{align}
	The third-order conditions for a general method cannot be solved using Mathematica. We therefore look for additional simplifications. If we impose $b_{i,i-1} = a_{i,i,i-1}$ (stiffly accurate) then the only two possible methods are:
	\begin{itemize}
		\item {\em method with rational coefficients}
		\begin{align}
			\begin{aligned}
				a_{221}&=1, &
				a_{321}&=-\frac{2}{3}, &
				a_{332}&=\frac{2}{3}, &
				a_{421}&=\frac{5}{12}, &
				a_{432}&=-\frac{5}{12}, \\
				a_{443}&=\frac{1}{2}, &
				a_{521}&=-\frac{1}{2}, &
				a_{532}&=\frac{1}{6}, &
				a_{543}&=\frac{2}{3}, &
				a_{554}&=\frac{2}{3}.
			\end{aligned}
			\label{eq:third-order-stiffly-accurate-a}
		\end{align}
		\item {\em method with non-closed form coefficients (shown in 32 digits of precision)}
		\begin{align}
			\begin{scriptsize}
			\begin{aligned}
				a_{110}&= 0.28752089074127389355949327393422 &
				a_{210}&= 0.28752089074127389355949327393422 \\
				a_{221}&= 4.5591287787686468000885595638962 &
				a_{310}&= 1.2842410657027653598348788927797 \\
				a_{321}&= 0.47597789152576385936348211385602 &
				a_{332}&= 0.092460512109263210293986365244118 \\
				a_{410}&= 0.66264895274771687678086938284766 &
				a_{421}&= 0.0017995404978923055219413775561107 \\
				a_{432}&= -0.040810958288454482287992176299806 &
				a_{443}&= 0.37636246504284529998518141589603.
			\end{aligned}
			\end{scriptsize}\ignorespaces
			\label{eq:third-order-stiffly-accurate-b}
		\end{align}
	\end{itemize}
	Alternatively if we impose singly-implicit structure ($a_{i,i,i-1} = \gamma$) then there exists a family of methods with arbitrary $\gamma$, however we are unable to obtain a close-form equations using Mathematica. One example method is \cref{eq:IMEX-NPRK354-Si}.

\section{IMIM-NPRK method derivations -- second-order}
\label{app:method-derivations-imim}
We start from the ansatz
\begin{align}
	\begin{aligned}
		Y_1 &= y_n, \\
		Y_2 &= y_n + h a_{221} F(Y_2,Y_1), \\
		Y_3 &= y_n + h a_{321} F(Y_2,Y_1) + h a_{323} F(Y_2,Y_3), \\
		y_{n+1} &= y_n + b_{21} F(Y_2,Y_1) + b_{23} F(Y_2,Y_3).
	\end{aligned}
\end{align}
The underlying method $M^{\{1\}}$ does not use the stage $Y_3$ while method $M^{\{2\}}$ does not use the stage $Y_2$. The reduced Tableaux for methods $M^{\{1\}}$ and $M^{\{2\}}$ are
\begin{center}
	\begin{tabular}{l|lll}
		 & $a_{221}$  \\ \hline
		 & $(b_{21} + b_{23})$
	\end{tabular}
	\hspace{5em}
	\begin{tabular}{l|ll}
		 & 0 \\
		 & $a_{321}$ & $a_{323}$ \\ \hline
		 & $b_{21}$ & $b_{23}$
	\end{tabular}	
\end{center}	
	and the order conditions for the NPRK method are
	\begin{align}
		b_{21} + b_{23} &= 1, &
		a_{221} (b_{21} + b_{23}) &= 1/2, &
		a_{321} b_{21} + a_{323} b_{23} = 1/2.
	\end{align}
	This leads to the two-parameter family of methods
	\begin{align}
		b_{21}  &= 1 - b_{23}, &
		a_{221} &= \frac{1}{2}, &
		a_{321} &= \frac{1 - 2a_{323}b_{23}}{2 - 2b_{23}},
	\end{align}
	with linear stability function
	\begin{align*}
		R(z_1,z_2) &= \frac{a_{221} z_1 (a_{323} z_2-1)+(z_1+z_2) (a_{321} b_{23} z_2-a_{323} b_{21} z_2+b_{21}+b_{23})-a_{323} z_2+1}{(a_{221} z_1-1) (a_{323} z_2-1)}.
	\end{align*}
	The infinite $z_1$ and $z_2$ limits are
	\begin{align*}
		\begin{aligned}
			\lim_{z_1 \to \infty} \alpha(z_1,z_2) &= \frac{a_{221} a_{323} z_2-a_{221}+a_{321} b_{23} z_2-a_{323} b_{21} z_2+b_{21}+b_{23}}{a_{221} a_{323} z_2-a_{221}}, \\
			\lim_{z_2 \to \infty} \alpha(z_1,z_2) &= \frac{\infty  (a_{321} b_{23}-a_{323} b_{21})}{a_{221} a_{323} z_1-a_{323}}.
		\end{aligned}
	\end{align*}
	To ensure a bounded amplification factor in the $z_2 \to \infty$ limit, we require the additional condition
	\begin{align}
		a_{321} b_{23}-a_{323} b_{21} = 0.
	\end{align}
	This leads to a one parameter family of methods
	\begin{align}
		b_{21}  &= 1-b_{23}, &
		a_{221} &= \frac{1}{2}, \\
		a_{321} &= \frac{1-b_{23}}{4 (b_{23}-1) b_{23}+2}, &
		a_{323} &= \frac{b_{23}}{4 (b_{23}-1) b_{23}+2},
	\end{align}
	with stability function 
	\begin{align}
		-\frac{z_2 (b_{23} (8 b_{23}+z_1-10)+4)+2 (2 (b_{23}-1) b_{23}+1) (z_1+2)}{(z_1-2) \left(4 b_{23}^2-b_{23} (z_2+4)+2\right)}.
	\end{align}
	We note that this is a rational function of the form
	\begin{align}
		\frac{(z_1 + c_1)(z_2 + c_2) + c_3}{(z_1-2)(z_2+d_0)},
	\end{align}
	where the constants depend on $b_{23}$ and are given by
	\begin{align}
		d_0 &= 4-4 b_{23}-\frac{2}{b_{23}}, \quad
		c_1 = 8 b_{23}+\frac{4}{b_{23}}-10, \quad
		c_2 =4 b_{23}+\frac{2}{b_{23}}-4, \\		
		c_3 &= -32 b_{23}^2-\frac{8}{b_{23}^2}+80 b_{23}+\frac{40}{b_{23}}-80.
	\end{align}
	If we force $c_3=0$ (this requires $b_{23} \in \left\{ \tfrac{1}{2}, 1 \right\}$) then the stability function is separable, in the sense that it can be written as 
	\begin{align}
		R(z_1,z_2) = f(z_1)g(z_2).
	\end{align}
	Moreover, for both choices of $b_{23}$ we have 
	\begin{align}
		\begin{aligned}
			f_1(z_1) = \frac{z_1+2}{z_1-2}, &&
			g_1(z_2) = \frac{z_2+2}{z_2-2},
		\end{aligned}
	\end{align}
	which implies that the resulting methods are A-stable. The corresponding methods are

	\begin{itemize}
		\item Implicit-Implicit Midpoint ($b_{23} = 1$) and its transpose (equivalent if $F(u,v) = F(v,u)$)
		\begin{align}
			\begin{aligned}
				Y_1 &= y_n \\
				Y_2 &= y_n + \frac{h}{2} F(Y_2,Y_1) \\
				Y_3 &= y_n + \frac{h}{2} F(Y_2,Y_3) \\
				y_{n+1} &= y_n + h F(Y_2,Y_3)
			\end{aligned}
			&&
			\begin{aligned}
				Y_1 &= y_n \\
				Y_2 &= y_n + \frac{h}{2} F(Y_1,Y_2) \\
				Y_3 &= y_n + \frac{h}{2} F(Y_3,Y_2) \\
				y_{n+1} &= y_n + h F(Y_3,Y_2)
			\end{aligned}
			\label{eq:imim-midpoint-with-transpose}	
		\end{align}
		$M^{\{1\}}$ and $M^{\{2\}}$ are both the implicit midpoint method
		\begin{center}	
			\renewcommand*{\arraystretch}{1.5}
			\begin{tabular}{l|lll}
				 & 0 \\
				 & 0 & $\tfrac{1}{2}$  \\ \hline
				 & 0 & $1$
			\end{tabular}
		\end{center}
		
		\item Implicit-Implicit Midpoint/Crank Nicolson ($b_{23} = \tfrac{1}{2}$) and its transpose
		\begin{align}
			\begin{small}
				\begin{aligned}
					Y_1 &= y_n \\
					Y_2 &= y_n + \frac{h}{2} F(Y_2,Y_1) \\
					Y_3 &= y_n + \frac{h}{2} F(Y_2,Y_1) + \frac{h}{2} F(Y_2,Y_3) \\
					y_{n+1} &= y_n + \frac{h}{2} F(Y_2,Y_1) + \frac{h}{2} F(Y_2,Y_3)
				\end{aligned}
			\end{small}
			&&
			\begin{small}
				\begin{aligned}
					Y_1 &= y_n \\
					Y_2 &= y_n + \frac{h}{2} F(Y_1,Y_2) \\
					Y_3 &= y_n + \frac{h}{2} F(Y_1,Y_2) + \frac{h}{2} F(Y_3,Y_2) \\
					y_{n+1} &= y_n + \frac{h}{2} F(Y_1,Y_2) + \frac{h}{2} F(Y_3,Y_2)
				\end{aligned}
			\end{small}
			\label{eq:imim-crank-with-transpose}		
		\end{align}
		$M^{\{1\}}$ is the implicit midpoint while $M^{\{2\}}$ is the implicit Crank Nicolson method
		\begin{center}	
			\renewcommand*{\arraystretch}{1.5}
			\begin{tabular}{l|lll}
				 & 0 \\
				 & 0 & $\tfrac{1}{2}$  \\ \hline
				 & 0 & $1$
			\end{tabular}
			\hspace{5em}	
			\begin{tabular}{l|lll}
				 & 0 \\
				 & $\tfrac{1}{2}$ & $\tfrac{1}{2}$  \\ \hline
				 & $\tfrac{1}{2}$ & $\tfrac{1}{2}$
			\end{tabular}
		\end{center}	
	\end{itemize}
	
	\begin{remark}
		Non-separable solutions with $b_{23} > 1/2$ appear to also have good stability properties.	
	\end{remark}
	
	\vfill
	\begin{center}
		\includegraphics[width=1em]{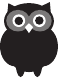}	
	\end{center}